\documentclass[11pt,a4paper]{article}
\usepackage[utf8]{inputenc}
\usepackage[english]{babel}
\usepackage[T1]{fontenc}
\usepackage{amscd,amssymb,amsmath,amsthm,amstext,amsfonts}
\usepackage{stmaryrd}
\usepackage[pdftex]{graphicx}
\usepackage{fancybox}
\usepackage[square,numbers]{natbib}
\usepackage{booktabs}
\usepackage{hyperref}
\usepackage{algorithm}
\usepackage{algpseudocode}
\usepackage{enumerate}  
\usepackage[capitalise,nameinlink]{cleveref}
\crefformat{equation}{\textup{#2(#1)#3}}
\crefname{figure}{Figure}{Figures}
\crefformat{subfigure}{Figure #2#1#3}
\usepackage{comment}
\usepackage{mathtools}
\usepackage{tikz,pgfplots}
\pgfplotsset{compat=newest}
\usepackage{pgfplotstable}
\usepackage{caption}
\usepackage{subcaption}
\pgfplotsset{compat=newest}
\usepackage{multirow}
\usepackage{authblk} 
\usepackage[normalem]{ulem}
\usepackage{amsopn}

\newcommand{\M}{\mathcal{M}}
\newcommand{\RE}{\mathbb{R}}

\newtheorem{theorem}{Theorem}[section]

\theoremstyle{definition}
\newtheorem{definition}[theorem]{Definition}
\theoremstyle{remark}
\newtheorem{remark}[theorem]{Remark}
\numberwithin{theorem}{section}
\numberwithin{equation}{section}
\numberwithin{table}{section}
\numberwithin{figure}{section}

\allowdisplaybreaks

\author[1]{Moataz M. Alghamdi}
\affil[1]{Applied Mathematics and Computational Sciences (AMCS), King Abdullah University of Science and Technology. Thuwal, 23955-6900, Kingdom of Saudi Arabia}

\author[1,2]{Daniele Boffi}
\affil[2]{Department of Mathematics ``F. Casorati'', University of Pavia. Via Ferrata 1, 27100 Pavia, Italy}

\author[3]{Francesca Bonizzoni}
\affil[3]{MOX - Department of Mathematics, Politecnico di Milano. Via Bonardi 9, 20133 Milano, Italy}

\title{A greedy MOR method for the tracking of eigensolutions to parametric elliptic PDEs}

\date{}

\begin{document}

\maketitle

\begin{abstract}
In this paper we introduce an algorithm based on a sparse grid adaptive refinement, for the approximation of the eigensolutions to parametric problems arising from elliptic partial differential equations.
In particular, we are interested in detecting the crossing of the hypersurfaces describing the eigenvalues as a function of the parameters.

The a priori matching is followed by an a posteriori verification, driven by a suitably defined error indicator. At a given refinement level, a sparse grid approach is adopted for the construction of the grid of the next level, by using the marking given by the a posteriori indicator.

Various numerical tests confirm the good performance of the scheme.
\end{abstract}

\vspace{1cm}


%%%%%%%%%%%%%%%%%%%%%%%%%%%%%%%%%%%%%%%%%%%%%%%%%%%%%%%%%%%%%%%%%%%%%%%%%%%%%%%%%%
\section{Introduction} \label{se:Introduction}
%%%%%%%%%%%%%%%%%%%%%%%%%%%%%%%%%%%%%%%%%%%%%%%%%%%%%%%%%%%%%%%%%%%%%%%%%%%%%%%%%%

Many engineering applications require the knowledge of resonance frequencies of the considered structure. Prime examples are vibration problems in mechanical engineering, where the vibration of buildings or bridges at their natural frequencies might cause damage and structural failure.

During the last decades, computation power has substantially increased. Nevertheless, the computation effort required for the solution of large-scale eigenproblems is still considerable. 
The situation gets more difficult when manufacturing imperfections and geometric or material variability are included in the mathematical model as parameters or random fields. In particular, when the solution of the same eigenproblem is of interest for many different values of the parameters, the direct computation of eigenpairs by means of standard numerical techniques entails an unaffordable computational effort.

Model Order Reduction (MOR) methods aim at reducing the overall computational effort by producing a surrogate of the parameter-to-eigenpair map, which is accurate and at the same time fast to be evaluated. They develop in two phases: several snapshots of eigenpairs (corresponding to an appropriate set of parameter values) are computed during the offline phase, and used to construct the surrogate, which is then evaluated at any new parameter value during the online phase.

A very popular ROM method developed during the last few decades is the Reduced Basis (RB) method (see, e.g., \cite{Quarteroni-Manzoni-Negri-2015, Hesthaven-Rozza-Stamm-2016}) which constructs the surrogate by projection onto the set of precomputed snapshots selected either adaptively by means a greedy algorithm or by Proper Orthogonal Decomposition (POD). 

Up to now, a tremendous effort has been made in developing ROM for source problems, whereas model reduction for parametric/stochastic eigenproblems is still a largely unexplored field. The way was paved by the pioneering work~\cite{Machiels-Maday-Oliveira-Patera-Rovas-2000}, where an RB approach is proposed to approximate the smallest eigenvalue. This methodology has been further developed to deal with more eigenvalues in~\cite{Pau-2007}. A component-based RB method for eigenvalue problems is proposed in~\cite{Vallaghe-Huynh-Knezevic-Nguyen-Patera-2015} and an a posteriori estimator for eigenvalues is studied in~\cite{Huynh-Knezevic-Patera-2013}. All the previously mentioned contributions do not cover the case of multiple eigenpairs. Instead, in~\cite{Horger-Wohlmuth-Dickopf-2017}, an a posteriori error bound for multiple eigenvalues (but not eigenvectors) is studied under the assumption of affine parametric dependence of the eigenproblem. Finally, in~\cite{Fumagalli-Manzoni-Prolini-Verani-2016} a greedy RB method for both affine and non-affine parametric eigenproblems is proposed, with a focus on the smallest (single) eigenpair, only. The same task (namely, the approximation of the first eigenpair) in the context of stochastic eigenproblems has been solved in~\cite{Hakula-Kaarnioja-Laaksonen-2015} and~\cite{Andreev-Schwab-2012} by means of the Stochastic Galerkin and Stochastic Collocation method, respectively. 

The aim of the present paper is the development of an algorithm to match (and interpolate) snapshots of eigenpairs of symmetric problems as the parameter varies in a given $p$-di\-men\-sion\-al subset of interest $\mathcal M\subset\RE^p$, and as the eigenvalues lie in a fixed window of interest $I_\lambda\subset\RE$. In particular, the issues connected with multiple eigenvalues as well as crossings of the hypersurfaces described by the eigenvalues as the parameter varies in $\mathcal M$ are thoroughly analyzed. 

The parameter space is sampled adaptively, following a two-phase procedure. First, a suitably adapted version of the a priori matching proposed in~\cite{Nobile-Pradovera-2021} is applied. Numerical examples (see Section~\ref{se:BuildingBlock}) show that this technique might produce an incorrect matching. Hence, we have developed a novel a posteriori indicator based on the orthogonality of the snapshots of eigensolutions, which drives the adaptive sampling. The algorithm is first presented on a one-dimensional (in the parametric space) case, and then extended to the high-dimensional setting by means of hierarchical locally refined sparse grids.

It is worth mentioning that the a priori matching that we use is connected with MOR techniques developed in a different framework, namely, the parametric-in-frequency Helmholtz boundary value problem. In particular, we mention~\cite{BonizzoniNobilePerugia18,BonizzoniNobilePerugiaPradovera20a,BonizzoniNobilePerugiaPradovera20b,BonizzoniPradoveraRuggeri}, where rational-based surrogates for the Helmholtz solution map are constructed. Indeed, the roots of the denominator of the surrogate are approximations to the resonances of the Helmholtz problem, which in turn are eigenvalues of the corresponding elliptic problem. 
A matching strategy in the same spirit as~\cite{Nobile-Pradovera-2021} and the a priori matching of the present paper is studied in~\cite{Yue-Feng-Benner-2019}. Moreover, a mode tracking method for the parametrized Maxwell eigenvalue problem is presented in~\cite{GeorgAckermannCornoSchops, ZieglerGeorgAckermannSchoeps,YangAjjarapu,JorkowskiSchuhmann}.

Our work is the first step toward a deep understanding of Uncertainty Quantification for stochastic eigenvalue problems arising from elliptic PDEs with stochastic diffusion coefficients. This field of research, even though relatively unexplored, is extremely important in various application areas.

The paper is organized as follows. In Section~\ref{se:Setting} we describe the problem of interest. Section~\ref{se:Algorithm} is dedicated to the adaptive algorithm: first we give a general overview of the main steps of the proposed algorithm, details on all the steps then follow. In Section~\ref{se:numerical_results} we present both one and two dimensional results to validate the proposed strategy. 
Conclusions are finally drawn in Section~\ref{se:conclusions}.

%%%%%%%%%%%%%%%%%%%%%%%%%%%%%%%%%%%%%%%%%%%%%%%%%%%%%%%%%%%%%%%%%%%%%%%%%%%%%%%%%%
\section{Setting of the Problem} \label{se:Setting}

Let us consider the Hilbert triplet
\[
V\subset H\simeq H'\subset V'
\]
and a parameter space $\M\subset\RE^d$.
It is out of the aims of this work to identify the most general assumptions on the parameter space. Very often in the applications we have in mind, it is a tensor product of intervals (like a hypercube). In any case in what follows we will need that it is connected and that it supports an initial grid with neighboring points as described later in this section.

For each $\mu\in\M$, we consider two symmetric and bilinear forms
\[
\aligned
&a(\cdot,\cdot;\mu):V\times V\to\RE,\\
&b(\cdot,\cdot;\mu):H\times H\to\RE.
\endaligned
\]
Furthermore, we make the following assumptions
\[
\aligned
&V\text{ compact in }H,\\
&a(\cdot,\cdot;\mu)\text{ elliptic in }V\quad\forall\mu\in\M,\\
&b(\cdot,\cdot;\mu)\text{ equivalent to the inner product in }H\quad\forall\mu\in\M.\\
\endaligned
\]
More general assumptions could be made. However, the features of our strategy are better described in this simpler setting.

Our aim is to approximate the solutions of the following parametric eigenvalue problem: for all $\mu\in\M$, find real eigenvalues $\lambda(\mu)$ and non-vanishing eigenfunctions $u(\mu)\in V$ such that
\begin{equation}  
\label{eq:main}
a(u,v;\mu)=\lambda(\mu)b(u,v;\mu)\quad\forall v\in V.
\end{equation}

Our assumptions ensure that the problem is associated with a compact solution operator so that all eigenvalues $\{\lambda_j(\mu)\}_{j=1}^\infty$ correspond to finite-dimensional eigenspaces.
In particular, we are interested in detecting the behavior of the hypersurfaces defined by $\lambda_j(\mu)$ in the region $\M\times\RE$. These hypersurfaces may intersect, leading in general to multiple eigenvalues at one point of intersection. Indeed, the situation can be very complicated when the dimension $d$ of the parametric space gets large. In the simplest case, $d=1$, we are dealing with intersections of curves.

To better define our problem, we restrict the range of the eigenvalues we are interested in to an interval $I_\lambda=[\lambda_{\min},\lambda_{\max}]$, also referred to as \emph{window of interest}. 
Consequently, we only examine the hypersurfaces in the region $\M\times I_\lambda$. This implies in particular that hypersurfaces can enter or exit this region of interest when the corresponding eigenvalues cross the values of $\lambda_{\min}$ or $\lambda_{\max}$. It follows then that the number of eigenvalues considered for $\mu_1$ and $\mu_2$ in $\M$ can be different from each other when $\mu_1\ne\mu_2$.

Problem~\eqref{eq:main} is discretized by finite elements (FEs). That is, we consider a finite dimensional subspace $V_h\subset V$ and for all $\mu\in\M$ we consider the matrix generalized eigenvalue problem: find real eigenvalues $\lambda_h(\mu)$ and non-vanishing eigenfunctions $u_h(\mu)\in V_h$ such that
\begin{equation}
    \label{eq:mainh}
    a(u_h,v;\mu)=\lambda_h(\mu)b(u_h,v;\mu)\quad\forall v\in V_h.
\end{equation}
This will be considered as our \emph{high-fidelity} solution.

%%%%%%%%%%%%%%%%%%%%%%%%%%%%%%%%%%%%%%%%%%%%%%%%%%%%%%%%%%%%%%%%%%%%%%%%%%%%%%%%%%
\section{Description of the adaptive algorithm} 
\label{se:Algorithm}
%%%%%%%%%%%%%%%%%%%%%%%%%%%%%%%%%%%%%%%%%%%%%%%%%%%%%%%%%%%%%%%%%%%%%%%%%%%%%%%%%%

Before examining the details of our approach, we give a general overview of the main steps taken to track the behavior of the hypersurfaces in the region $\M\times I_\lambda$ described by the varying eigenvalues. We assume that we are given a grid in the parametric space $\M$ and that we have computed the eigenvalues in the interval $I_\lambda$ as well as the corresponding eigenfunctions for each point of the grid.

The first phase consists in the \textbf{a priori matching} of the eigenvalues of each pair of neighboring parameters $\mu_i$ and $\mu_k$ where the indices $i$ and $k$ are defined based on the used grid. The matching is performed by considering, in a suitable sense, how close the eigenvalues and the eigenfunctions are to each other.
This phase is prone to error in particular if the distance between $\mu_i$ and $\mu_k$ is large with respect to the variability of the eigenvalues.
The a priori phase is followed by an \textbf{a posteriori verification} of the matching. We introduce a suitable a posteriori indicator that is based on the orthogonality of the eigenfunctions. This phase aims to confirm whether the a priori matching was performed correctly or not. If not, the corresponding interval is marked for refinement.
Finally, a sparse grid approach is used to drive the \textbf{refinement strategy}, leading to an adaptive procedure that is terminated when a suitable stopping criterion is met.

The three phases of our adaptive strategy are executed using three interrelated algorithms. Algorithm~\ref{alg:refinement} describes the global refinement procedure. In turn, this algorithm relies on two additional algorithms. The first one, Algorithm~\ref{alg:a-priori}, describes the local a priori matching procedure while the second one, Algorithm~\ref{alg:a-posteriori2}, enforces the local a posteriori test.

%%%%%%%%%%%%%%%%%%%%%%%%%%%%%%%%%%%%%%%%%%%%%%%%%%%%%%%%%%%%%%%%%%%%%%%%%%%%%%%%%%
\subsection{The a Priori Matching}
%%%%%%%%%%%%%%%%%%%%%%%%%%%%%%%%%%%%%%%%%%%%%%%%%%%%%%%%%%%%%%%%%%%%%%%%%%%%%%%%%%

The FE method~\eqref{eq:mainh} leads to an algebraic problem as follows: find $\lambda(\mu) \in \mathbb{R}$ and $ u(\mu) \in \mathbb{R}^N$ with $u(\mu) \neq 0 $ such that $A(\mu)u(\mu) = \lambda(\mu) B(\mu)u(\mu)$, where $A(\mu)$ and $B(\mu)$ are matrices in $\mathbb{R}^{N \times N}$ where $N$ is the dimension of our finite element space. $A(\mu)$ and $B(\mu)$ are symmetric and positive definite for all values of $\mu \in \M$. To avoid heavy notation, we denote discrete quantities without indicating the space mesh index $h$ ($\lambda$ instead of $\lambda_h$, etc.).

We consider two different parameter values $\mu_i$ and $\mu_k$ that are the endpoints of what we are going to call from now on a \emph{local subinterval}. For these two points, we generate the two sets of FE eigenpairs $\{(\lambda_j(\mu_i),u_j(\mu_i))\}_{j=1}^{n_i}$, $\{(\lambda_\ell(\mu_k),u_\ell(\mu_k))\}_{\ell=1}^{n_k}$. Note that the values $n_i$ and $n_k$ may be different from each other in particular since we are looking for all eigenvalues within the window of interest $I_\lambda$. We assumed that the bilinear form $b(\cdot,\cdot,\mu)$ is equivalent to the scalar product in $H$ for all $\mu\in\M$; we denote the associated norm by
\[
\|v\|_{b,\mu}=b(v,v,\mu)^{1/2}
\]
The eigenfunctions are normalized with respect to this norm so that
\[
\|u_j(\mu_i)\|_{b,\mu_i}=\|u_\ell(\mu_k)\|_{b,\mu_k}=1.
\]

\begin{definition}
\label{def:cost_matrix}
The cost matrix $D^{(i,k)}$ associated with the local subinterval with endpoints $\mu_i$ and $\mu_k$ has size $n_i\times n_k$ and entries
\begin{equation}
    \label{eq:cost_matrix}
    D^{(i,k)}_{j,\ell}\coloneqq w_1 |\lambda_j(\mu_i) - \lambda_\ell(\mu_k) |
        + w_2 \min\{ \|u_j(\mu_i)-u_\ell(\mu_k)\|_{b,\bar\mu} , \|u_j(\mu_i)+u_\ell(\mu_k)\|_{b,\bar\mu}\},
\end{equation}
where $w_1,w_2\in\mathbb R_+$ are weights and $\bar\mu$ is a fix parameter value between $\mu_i$ and $\mu_k$.
\end{definition}

\begin{remark}
A similar definition to~\eqref{eq:cost_matrix} can be found in~\cite{Nobile-Pradovera-2021}. However, some modifications are necessary in the present context, since the normalization of the eigenfunctions doesn't necessarily imply the same choice of sign. Hence we have to compare both the sum and the difference of corresponding eigenfunctions.
\end{remark}

Once we have introduced the cost matrix, we want to minimize its entries by solving the following optimization problem: find the permutation $\sigma^\star=(\sigma_1,\ldots,\sigma_{\bar n})\in(1,\ldots,\bar n)!$ such that
\begin{equation}
    \label{eq:min_cost_matrix}
    \sigma^\star\coloneqq\operatorname{argmin}_{\sigma\in(1,\ldots,\bar n)!}
    \sum_{\alpha=1}^{n_i} \sum_{\beta=1}^{n_k} D^{i,k}_{\sigma_\alpha,\sigma_\beta},
\end{equation}
with $\bar n\coloneqq \min\{n_i,n_k\}$.
Between the several options available to compute the solution of~\eqref{eq:min_cost_matrix}, we adopt the the so-called Hungarian Algorithm~\cite{Kuhn1955}. The permutation solution is then used to reorder the eigensolutions so that there is a one-to-one correspondence between the matched eigensolutions

The structure of our code is reported in the following Algorithm~\ref{alg:a-priori}, where we assume, without loss of generality, that $n_i\geq n_k$, i.e., $\bar n= n_k$.

\begin{algorithm}
\caption{Local a priori matching}
\label{alg:a-priori}
\begin{algorithmic}[1]
    \Require{$\mu_i,\,\mu_k\in\mathcal M$,  
        $\{(\lambda_j(\mu_i),u_j(\mu_i))\}_{j=1}^{n_i}$, 
        $\{(\lambda_\ell(\mu_k),u_\ell(\mu_k))\}_{\ell=1}^{n_k}$, $w_1,w_2\in\mathbb R_+$}
    \Ensure{Reordered eigenpairs $\{(\lambda^\star_j(\mu_i),u^\star(\mu_i))\}_{j=1}^{n_i}$, 
        $\{(\lambda^\star_\ell(\mu_k),u^\star_\ell(\mu_k))\}_{\ell=1}^{n_k}$ }
    \State{Compute $D^{(i,k)}\in\mathbb R^{n_i\times n_k}$ with weights $w_1,\, w_2$}
        \Comment{See~\eqref{eq:cost_matrix}}
    \State{Find $\boldsymbol{\sigma}^\star$ solution to~\eqref{eq:min_cost_matrix}}
        \Comment{Hungarian algorithm}   
    \State{Set $(\lambda^\star_j(\mu_i),u^\star(\mu_i))=(\lambda_j(\mu_i),u_j(\mu_i))$ for $j=1,\ldots,n_i$} \State{Set $(\lambda^\star_\ell(\mu_k),u^\star_\ell(\mu_k))=(\lambda_{\sigma^\star_\ell}(\mu_k),u_{\sigma^\star_\ell}(\mu_k))$ for $\ell=1,\ldots,n_k$}
\end{algorithmic}
\end{algorithm}

%%%%%%%%%%%%%%%%%%%%%%%%%%%%%%%%%%%%%%%%%%%%%%%%%%%%%%%%%%%%%%%%%%%%%%%%%%%%%%%%%%
\subsection{The a Posteriori Verification}
\label{sec:a_poteriori_verification}
%%%%%%%%%%%%%%%%%%%%%%%%%%%%%%%%%%%%%%%%%%%%%%%%%%%%%%%%%%%%%%%%%%%%%%%%%%%%%%%%%%

After the a priori matching, we introduce an a posteriori verification phase, which is based on the following projection matrix. 

\begin{definition}
The projection matrix $\Pi^{(i,k)}$ associated with the local subinterval with endpoints $\mu_i$ and $\mu_k$ has size $n_i\times n_k$ and its
entries are given by
\begin{equation}
    \label{eq:proj_matrix}
    \Pi^{(i,k)}_{j,\ell}=|b(u_j(\mu_i),u_\ell(\mu_k),\bar\mu)|,
\end{equation}
with $\bar\mu$ being a fixed parameter value between $\mu_i$ and $\mu_k$.
\end{definition}

Ideally, if the matching was performed correctly, the projection matrix should be close to diagonal: two matching eigenfunctions should be similar to each other and two non-matching eigenfunctions should be close to orthogonal with respect to the bilinear form $b$. In order to check whether the projection matrix is close to diagonal, we make use of a positive tolerance value $t_\pi$. This is used to truncate $\Pi^{(i,k)}$ --- inside a loop over $j=1,\ldots,\min\{n_i,n_k\}$ ---  in accordance with lines 4-13 of Algorithm~\ref{alg:a-posteriori2}. Let $r_1, r_2$ be the vectors containing the non-zero elements of the $j$-th column and $j$-th row of $\Pi^{(i,k)}$, respectively (line 14). If both $r_1$ and $r_2$ contain just one element, the a priori matching is considered correct. Instead, if their length differs, the interval is marked for refinement (lines 16-19). Checking the orthogonality of eigenfunctions is a good stopping criterion in general, but might fail when we are close to multiple eigensolutions. In such a case, the orthogonality between distinct eigenfunctions depends on the solver and is not immediate to check in practice. For this reason, we introduce a second positive tolerance value $t_\lambda$ that is responsible for verifying if two (or more) eigenvalues belong to a cluster. This scenario corresponds in Algorithm~\ref{alg:a-posteriori2} to the case where $r_1$ and $r_2$ have the same length, larger than $1$. Lines 20-25 introduce a specific definition for when multiple eigensolutions are to be considered as one cluster of indistinguishable eigenfunctions. 
The local subinterval is marked for refinement when it fails the $t_\lambda$ stopping criterion.

\begin{remark}[Choice of the tolerances $t_\pi, t_\lambda$]
Note that in the limit $t_{\pi} \to 0$, the projection matrix $\Pi^{(i,k)}$ is diagonal. 
On the other hand, the truncation process will leave the projection matrix unchanged as $t_{\pi}  \to 1$, possibly leading to an infinite loop over the refinement level (see Section~\ref{sec:refinement}). 
There is then a threshold between between $t_{\pi}$ being large enough to capture potential errors in matching choices, and small enough to minimize the number of refinements. The selection of the optimal value for $t_\pi$ becomes more delicate as the dimension of the parameter space gets larger.

A similar balance must hold for $t_\lambda$. When there is considerable overlap between two (or more) eigenfunctions, this will translate into non-zero off-diagonal elements of $\Pi^{(i,k)}$. This can lead to a very large number of refinements in order for the subinterval to be certified, unless $t_\lambda$ is chosen such that those overlapping eigenfunctions are identified as a cluster. However, if $t_\lambda$ is chosen to be too large, then wrong matching choices of orthogonal eigenfunctions will be certified by the a posteriori estimator because such eigenfunctions will be incorrectly considered indistinguishable. This premature termination can lead to wrong results.
\end{remark}

\begin{algorithm}
\caption{Local a posteriori verification}
\label{alg:a-posteriori2}
\begin{algorithmic}[1]
    \Require{$\mu_i,\,\mu_k\in \mathcal M$,  
        $\{(\lambda_j(\mu_i),u_j(\mu_i))\}_{j=1}^{n_i}$, 
        $\{(\lambda_\ell(\mu_k),u_\ell(\mu_k))\}_{\ell=1}^{n_k}$, $t_\pi,\, t_\lambda\in\mathbb{R}_+$}
    \Ensure{$\operatorname{if\_refinement}=0$ or $\operatorname{if\_refinement}=1$} 
    \State{Set if\_refinement=0}
    \State{Compute $\Pi^{(i,k)}\in\mathbb R^{n_i\times n_k}$}
        \Comment{See Equation~\eqref{eq:proj_matrix}}
\For{$j=1$ to $\texttt{min}(n_i,n_k)$}
\For{$\ell=1$ to $n_k$}
    \If{$\Pi^{(i,k)}_{j,j}\ge\Pi^{(i,k)}_{j,\ell}+t_\pi$}
        \State Set $\Pi^{(i,k)}_{j,\ell}=0$
        \Comment{Truncate $\Pi^{(i,k)}_{j,:}$ up to tolerance $t_\pi$}
    \EndIf
\EndFor
\For{$\ell=1$ to $n_i$}
    \If{$\Pi^{(i,k)}_{j,j}\ge\Pi^{(i,k)}_{\ell,j}+t_\pi$}
        \State Set $\Pi^{(i,k)}_{\ell,j}=0$
        \Comment{Truncate $\Pi^{(i,k)}_{:,j}$ up to tolerance $t_\pi$}
    \EndIf
\EndFor
\State{Let $r_1=\texttt{find}(\Pi^{(i,k)}_{j,:})$,
    $r_2=\texttt{find}(\Pi^{(i,k)}_{:,j})$}
    \Comment{Indices of the non-zero elements}
\If{$(\texttt{length}(r_1) )> 1$ or $(\texttt{length}(r_2) )> 1$) }
    \If{$\texttt{length}(r_1)\neq\texttt{length}(r_2)$}
        \State{$\operatorname{if\_refinement}=1$}
         \Comment{Mark $[\mu_i,\mu_k]$ for refinement}
        \State{\textbf{break}}
    \EndIf
    \State{Let $\alpha^i_{j,\gamma_1}\coloneqq \frac{|\lambda_j(\mu_i)-\lambda_{\gamma_1}(\mu_i)|}{\lambda_j(\mu_i)}$ for $\gamma_1\in r_1$}
    \State{Let $\alpha^k_{j,\gamma_2}\coloneqq \frac{|\lambda_j(\mu_k)-\lambda_{\gamma_2}(\mu_k)|}{\lambda_j(\mu_k)}$ for $\gamma_2\in r_2$}
    \If{$\texttt{max}_{\gamma_1\in r_1,\, \gamma_2\in r_2}\{ \alpha^i_{j,\gamma_1},\alpha^k_{j,\gamma_2}\}>t_\lambda$}
    \Comment{No cluster is identified}
        \State{$\operatorname{if\_refinement}=1$}
         \Comment{Mark $[\mu_i,\mu_k]$ for refinement}
        \State{\textbf{break}}
    \EndIf

\EndIf

\EndFor
\end{algorithmic}
\end{algorithm}

%%%%%%%%%%%%%%%%%%%%%%%%%%%%%%%%%%%%%%%%%%%%%%%%%%%%%%%%%%%%%%%%%%%%%%%%%%%%%%%%%%
\subsection{The Sparse Grid-based Adaptive Sampling} 
\label{se:SparseGrid}
%%%%%%%%%%%%%%%%%%%%%%%%%%%%%%%%%%%%%%%%%%%%%%%%%%%%%%%%%%%%%%%%%%%%%%%%%%%%%%%%%%

In the present contribution we aim at a refinement strategy that performs well also when the dimension $d$ of our problem is large. To mitigate the curse of dimensionality, it is essential to pay particular attention to the way the grid of the parameter space $\M$ is refined. One possibility is to use locally-refined sparse grids, in the spirit of what was proposed in~\cite{Alsayyari-Perko-Lathouwers-Kloosterman-2019,Nobile-Pradovera-2021}. In particular, the refinement step described in Section~\ref{sec:refinement} relies on some notions and features related to sparse grids. For the readers' convenience, we recall them here.

The main ideas of our sparse grid approach are better explained, without loss of generality, when $\M=[-1,1]^2$. With small modifications, the case $\mathcal M=[a,b]\times [c,d]$, for $a,b,c,d\in\RE$ can be handled. Moreover, the following discussion can be easily extended to the high-dimensional framework, i.e., for $d\geq 3$.

\begin{figure}
    \centering
    \includegraphics[width=4cm]{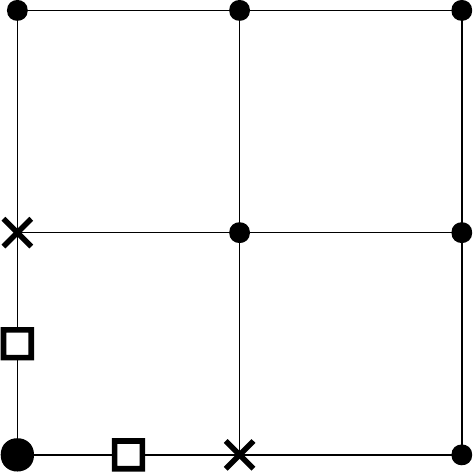}
    \includegraphics[width=4cm]{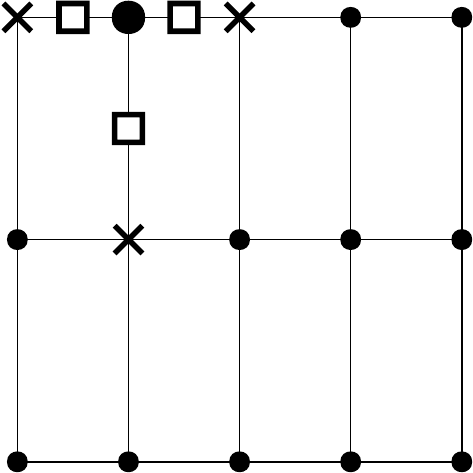}
    \includegraphics[width=4cm]{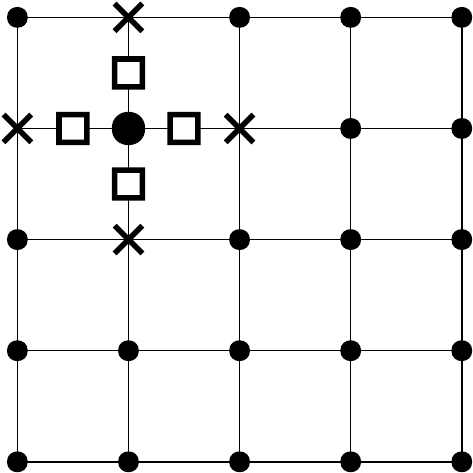}
    \caption{Forward points (squares) and neighbours (crosses) of $\mu=(-1,-1)\in\Gamma(1,1)$ (left); $\mu=(-\frac{1}{2},1)\in\Gamma(2,1)$ (middle); $\mu=(-\frac{1}{2},\frac{1}{2})\in\Gamma(2,2)$ (right).}
    \label{fig:forward_point}
\end{figure}

Let us define the sequence $\{\Gamma(m)\}_{m\in\mathbb N_0}$ of nested sets of points in $[-1,1]$ as follows:
\begin{equation}
    \Gamma(m)\coloneqq\left\{\begin{array}{ll}
        \{0\}&\text{if }m=0,\\
        \{2^{1-m }j\}_{j=-2^{m-1}}^{2^{m-1}}&\text{if }m>0.\\
    \end{array}\right.
\end{equation}
By tensor product, we get grids of points in $\mathcal M$. In particular, given a two-dimensional multi-index $\mathbf{m}=(m_1,m_2)\in\mathbb N_0^2$, the corresponding two-dimensional grid is defined as
\begin{equation*}
    \Gamma(\mathbf{m})=\Gamma(m_1)\times\Gamma(m_2)
    =\left\{(\alpha_1,\alpha_2),\, \alpha_k\in\Gamma(m_k),\, k=1,2\right\}.
\end{equation*}
Let $\mu=(\mu_1,\mu_2)\in\Gamma(\mathbf m)$ be given, and assume that both its entries are fractions in lowest terms. We define the set of \emph{forward points} $\mathcal V(\mu)$ of $\mu$ as 
\begin{equation}
    \label{eq:forward_points}
    \mathcal V(\mu)\coloneqq\{(\mu_1\pm 2^{-m_1},\mu_2),\, (\mu_1,\mu_2\pm 2^{-m_2})\}\cap \mathcal M,
\end{equation}
and the set of neighbours $\mathcal U(\mu)$ of $\mu$ as 
\begin{equation}
    \label{eq:neighbour}
    \mathcal U(\mu)\coloneqq\{(\mu_1\pm 2^{-(m_1-1)},\mu_2),\, (\mu_1,\mu_2\pm 2^{-(m_2-1)})\}\cap \mathcal M.
\end{equation}
Note that both sets $\mathcal V(\mu)$ and $\mathcal U(\mu)$ contain up to $2^d$ points for any $d\geq 2$. Some examples are depicted in Figure~\ref{fig:forward_point}. 

%%%%%%%%%%%%%%%%%%%%%%%%%%%%%%%%%%%%%%%%%%%%%%%%%%%%%%%%%%%%%%%%%%%%%%%%%%
\subsection{The Refinement Strategy}
\label{sec:refinement}
%%%%%%%%%%%%%%%%%%%%%%%%%%%%%%%%%%%%%%%%%%%%%%%%%%%%%%%%%%%%%%%%%%%%%%%%%%

We start with an initial partition $P^{(0)}$ or the parameter space $\M$. The procedure described in Algorithm~\ref{alg:refinement} describes how to go from the $\ell$-th partition $P^{(\ell)}$ to the next partition $P^{(\ell+1)}$ containing $P^{(\ell)}$.

Knowing that $P^{(\ell)}=P^{(\ell-1)}\cup P^{(\ell)}_\delta$, the refinement strategy takes into account only the points in $P^{(\ell)}_\delta$. This is essential in order to mitigate the curse of dimensionality when the value of $d$ is large. On the other hand, we do not need to consider the points in $P^{(\ell)}\setminus P^{(\ell)}_\delta$ since we know from the previous level that the a posteriori indicator confirmed the a priori matching there.

Each element (parameter) in $P^{(\ell)}_\delta$ has a collection of neighboring points (see equation~\eqref{eq:neighbour}). This defines a collection of local subintervals on which we apply Algorithms~\ref{alg:a-priori} and~\ref{alg:a-posteriori2}, leading to either a marking of the subinterval for further refinement or not. The output of one iteration of this algorithm defines the output of a \textbf{level}. We perform this algorithm iteratively for increasing values of the level $\ell\in\mathbb N_0$ (with the convention $P^{(-1)}=\emptyset$, so that $P^{(0)}=P^{(0)}_\delta$) until no extra refinements for any subinterval within the current level take place. 

\begin{algorithm}
\caption{Refinement Step}
\label{alg:refinement}
\begin{algorithmic}[1]
    \Require{Level $\ell\in\mathbb N_0$, initial grid $P^{(\ell)}=P^{(\ell-1)}\cup P^{(\ell)}_\delta$ with  $P^{(\ell-1)}=\{\mu^{(\ell-1)}_s\}_{s=1}^{N_{\ell-1}}$ and $P^{(\ell)}_\delta=\{\mu^{(\ell,\delta)}_s\}_{s=1}^{N_{\ell,\delta}}$,
    and eigenpairs $\Lambda^{(\ell)}_\delta\coloneqq\{(\lambda_j(\mu^{(\ell,\delta)}_s),u_j(\mu^{(\ell,\delta)}_s)\}_{j=1}^{n_s}$}
    \Ensure{Refined grid $P^{(\ell+1)}$ and eigenpairs $\Lambda^{(\ell+1)}_\delta$}
    \State{Set $P^{(\ell+1)}_\delta=\emptyset$, $\Lambda^{(\ell+1)}_\delta=\emptyset$}
    \For{$s=1:N_{\ell,\delta}$}
        \Comment{Loop on the points in $P^{(\ell)}_\delta$}
        \For{$\nu_{(r,s)}\in\mathcal V(\mu^{(\ell,\delta)}_s)$}
        \Comment{Loop on the neighbours of $\mu^{(\ell,\delta)}_s$ \eqref{eq:neighbour}}
            \State{$\Lambda^\star_s\coloneqq\{(\lambda_j(\mu_s^{(\ell,\delta)}),u_j(\mu_s^{(\ell,\delta)}))\}_{j=1}^{n_s}\gets
            \{(\lambda_j(\mu_s^{(\ell,\delta)}),u_j(\mu_s^{(\ell,\delta)}))\}_{j=1}^{n_s}$}
            \State{$\Lambda^\star_r\coloneqq\{(\lambda^\star_j(\nu_{(r,s)}),u^\star_j(\nu_{(r,s)}))\}_{j=1}^{n_s}\gets
            \{(\lambda_j(\nu_{(r,s)}),u_j(\nu_{(r,s)}))\}_{j=1}^{n_r}$}
             \Comment{See Algorithm~\ref{alg:a-priori}}
            \State{$\operatorname{if\_refinement}$ = a-posteriori check on $\Lambda^\star_s$, $\Lambda^\star_r$}
                \Comment{See Algorithm~\ref{alg:a-posteriori2}}
            \If{$\operatorname{if\_refinement}=1$}
                \State{$P^{(\ell+1)}_\delta\gets \bar \nu_{(r,s)}$} \Comment{$\bar\nu_{(r,s)}$ is the corresponding forward point of $\mu_s^{(\ell,\delta)}$}
                \State{Compute $\bar\Lambda\coloneqq\{(\lambda_{j}(\bar \nu_{(r,s)}),u_j(\bar \nu_{(r,s)}))\}_{j=1}^{n_r}$}
                \State{$\Lambda^{(\ell+1)}_\delta\gets\bar\Lambda$}
            \EndIf
        \EndFor
    \EndFor
    \State{Set $P^{(\ell+1)} = P^{(\ell)}\cup P^{(\ell+1)}_\delta$}
\end{algorithmic}
\end{algorithm}

%%%%%%%%%%%%%%%%%%%%%%%%%%%%%%%%%%%%%%%%%%%%%%%%%%%%%%%%%%%%%%%%%%%%%%%%%%%%%%%%%%
\section{Numerical results} 
\label{se:numerical_results}
%%%%%%%%%%%%%%%%%%%%%%%%%%%%%%%%%%%%%%%%%%%%%%%%%%%%%%%%%%%%%%%%%%%%%%%%%%%%%%%%%%

Section~\ref{se:BuildingBlock} is devoted to a numerical example illustrating the behavior of the two local procedures (a priori matching and a posteriori verification). The global refinement algorithm is discussed in one- and two- dimensional numerical examples in Sections~\ref{se:1D} and~\ref{se:2D}.

In all the presented numerical experiments we consider the following setting: let $\Omega=[0,1]^2$ be the physical domain and let $V=H^1(\Omega)$ and $H=L^2(\Omega)$ endowed with the usual inner products and norms. For all $\mu\in\mathcal M\subset\RE^d$ we look for eigenpairs $(\lambda(\mu),u(\mu))\in \RE^+\times H^1(\Omega)$ such that
\begin{equation}
    \label{eq:eigenvalue-pde}
    \left\{\begin{array}{ll}
         \nabla \cdot (c({\mu}) \nabla u({\mu})) = \lambda({\mu}) u({\mu})&  \text{ in }\Omega,\\
         u({\mu})=0& \text{ on }\partial\Omega,
    \end{array}\right.
\end{equation}
where $c(\mu)$ is a matrix with size $2\times 2$ and positive definite for all possible values of the parameter $\mu$.
Integrating by parts, we find the weak formulation of~\eqref{eq:eigenvalue-pde} of the form~\eqref{eq:main} with the bilinear forms
\begin{gather*}
    \begin{aligned}
    a(w,v;\mu)&\coloneqq\int_\Omega c(\mu) \nabla w\cdot\nabla v\, dx,  \\
    b(w,v;\mu)&\coloneqq\int_\Omega w v\, dx. 
    \end{aligned}
\end{gather*}
The FE method on a sufficiently refined mesh is employed to numerically approximate the eigenvalues and eigenfunctions at fixed values of the parameter $\mu\in\mathcal M$.
All the computations are performed in \textsc{Matlab} on a laptop with four cores (eight logical processors), 16 GB of RAM. Furthermore, we make use of the Partial Differential Equation toolbox~\cite{matlabPDE} and Sparse Grids Kit~\cite{piazzola.tamellini:SGK}.

%%%%%%%%%%%%%%%%%%%%%%%%%%%%%%%%%%%%%%%%%%%%%%%%%%%%%%%%%%%%%%%%%%%%%%%%%%
\subsection{Operations on a subinterval} 
\label{se:BuildingBlock}
%%%%%%%%%%%%%%%%%%%%%%%%%%%%%%%%%%%%%%%%%%%%%%%%%%%%%%%%%%%%%%%%%%%%%%%%%%

In the first numerical test we take $\mathcal M=[0.4,1]$,
\begin{equation}
\label{eq:c_1D}
    c({\mu}) = 
\begin{pmatrix}
{\mu}^{-2} & 1 \\
1 & 0.7^{-2}
\end{pmatrix},
\end{equation}
and $I_\lambda=[0,270]$.
We first detail the local steps of the proposed algorithm, namely, the a priori matching and the a posteriori verification, on the subinterval with endpoints $\mu_1=0.4$ and $\mu_2=0.7$.
The FE method at $\mu_1$ (respectively $\mu_2$) delivers four (respectively nine) eigenvalues in $I_\lambda$ (and corresponding eigenvectors), given by
\begin{gather*}
    \begin{aligned}
        \lambda(\mu_1)&=(80.8,\ 137.9,\ 230.6,\ 265.9)^\top,\\
        \lambda(\mu_2)&= (38.2,\ 81.1, \ 109.7,\ 129.4,\ 188.6,\ 189.9,\ 214.8,\ 260.9,\ 261.9)^\top.
    \end{aligned}
\end{gather*}

The $4 \times 9$ cost matrix (with weights  $w_1 = 1$, $w_2 = 200$) reads as follows:
$$
D = 
\begin{bmatrix*}[r]
\textcolor{red}{57.7} & 282.2 & 311.6 & 318.8 & 390.6 & 391.5 & 415.6 & 462.8 & 463.3\\
381.9 & \textcolor{red}{189.4} & 202.6 & 290.1 & 333.1 & 324.8 & 359.5 & 396.6 & 406.6\\
473.2 & 431.9 & 403.4 & 288.7 & \textcolor{red}{204.8} & 323.0 & 220.3 & 313.1 & 310.7\\
509.9 & 359.2 & 290.3 & 418.4 & 360.0 & 345.3 & 333.8 & \textcolor{red}{278.3} & 286.7
\end{bmatrix*}
$$
To minimize the total cost, namely, to solve the optimization problem~\eqref{eq:min_cost_matrix}, we employ the Hungarian algorithm, whose output is the permutation $\sigma^\star=(1,\ 2,\ 5,\ 8)$ (the entries of the cost matrix identified by $\sigma^\star$ are depicted in red). As a consequence, the vectors $\{\lambda_j(\mu_1),\ j=1,\ldots,4\}$ and $\{\lambda_\ell(\mu_2),\ \ell=\sigma_1^\star,\ldots,\sigma_4^\star)\}$ are matched as follows:
$$
\begin{array}{ccc}
     80.8  & \leftrightarrow & 38.2 \\
     137.9 & \leftrightarrow & 81.1 \\
     230.6 & \leftrightarrow & 188.6 \\
     265.9 & \leftrightarrow & 260.9
\end{array}
$$ 
and the vector $\{\lambda_\ell(\mu_2),\ \ell=1,\ldots,9\}$ is reordered accordingly:
\[
\lambda^\star(\mu_2)= (38.2,\ 81.1,\ 188.6,\ 260.9,\ 109.7,\ 129.4,\ 189.9,\ 214.8,\ 261.9)^\top.
\]

The a priori matching is then followed by the a posteriori verification. Here, we start with the $4\times 9$ projection matrix.
$$
\Pi = 
\begin{bmatrix}
    0.997  &  0.000  &  0.001 &  0.000 &  0.000 &   0.070  &  0.000 &  0.018   & 0.014 \\
    0.000  &  0.778 &  0.000 &  0.052 &  0.622 &   0.000  &  0.071  &  0.000  & 0.000 \\
    0.030 &  0.000 &  0.663 &  0.000 &  0.000  &  0.564 &  0.000  &  0.481 &   0.007 \\
   0.000 &   0.617 &  0.000  &  0.051  &  0.778  &    0.000 &   0.092  &  0.000  &  0.000
\end{bmatrix}
$$
We fix $t_\pi=0.21$, and we enter the loop over $j$ (line 3 in Algorithm~\ref{alg:a-posteriori2}. For $j=1$, the truncation conditions in line 5, 10 are fulfilled, hence all the off-diagonal entries of the first row and column of $\Pi$ are truncated to 0. For $j=2$, instead, the same process identifies two non-zero elements both in the second row and in the second column of $\Pi$:

$$
\Pi = 
\begin{bmatrix*}[r]
0.997 & 0 & 0 & 0 & 0 & 0 & 0 & 0 & 0\\
0 & 0.778 & 0 & 0 & 0.622 & 0 & 0 & 0 & 0\\
0 &  0 &  0.663 &  0.000 &  0.000  &  0.564 &  0.000  &  0.481 &  0.007\\
0 &   0.617 &  0.000  &  0.051  &  0.778  &    0.000 &   0.092  &  0.000  &  0.000
\end{bmatrix*}
$$
As a consequence, the first matching choice $80.8 \leftrightarrow 38.2$ is certified but not the second matching $137.9 \leftrightarrow 81.1$.
Here, the second eigenvalue at $\mu_1$, namely 137.9, is permitted to be matched to both the second (81.1 since $\Pi_{2,2} = 0.778$) and fifth (109.7 since $\Pi_{2,5} = 0.622$) eigenvalues at $\mu_2$. Furthermore, the second eigenvalue at $\mu_2$, 81.1, can be matched to the second and fourth (265.9 since $\Pi_{4,2} = 0.617$) eigenvalues at $\mu_1$.

We then check whether the second and fourth eigenvalues at $\mu_1$ form a cluster given the fixed tolerance $t_\lambda = 0.001$. A quick computation yields the following 
$$
\frac{|\lambda_2(\mu_1)-\lambda_4(\mu_1)|}{ \lambda_2(\mu_1)}>0.001
$$ 
and similar results hold at $\mu_2$ for the second and fifth eigenvalues. Thus, these eigenvalues are distinguishable from each other and no eigenvalue clusters are identified. Consequently, the a posteriori verification does not confirm the results of the a priori matching, and the subinterval is not certified and is marked for further refinement.

%%%%%%%%%%%%%%%%%%%%%%%%%%%%%%%%%%%%%%%%%%%%%%%%%%%%%%%%%%%%%%%
\subsection{Full 1D example} 
\label{se:1D}
%%%%%%%%%%%%%%%%%%%%%%%%%%%%%%%%%%%%%%%%%%%%%%%%%%%%%%%%%%%%%%%

The main purpose of this section is to show how refinement in a full 1D parametric domain takes place. We continue using the same set-up presented in Section~\ref{se:BuildingBlock}.

\begin{figure}[htpb]
\centering
\includegraphics[width=6cm]{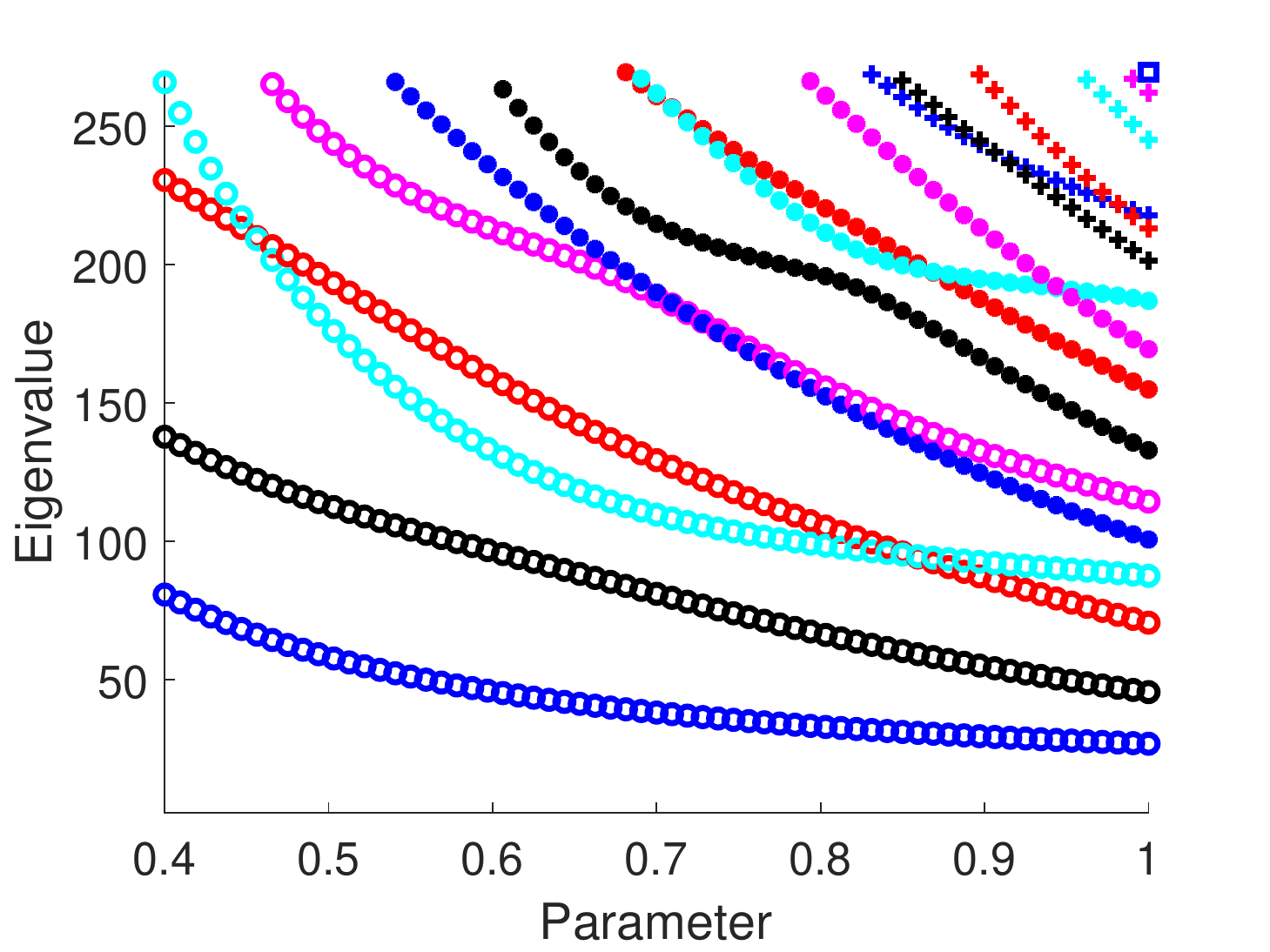}
\caption{Reference solution.} 
\label{fig:ReferenceSolution1D}
\end{figure}

For the sake of comparison, we numerically compute the reference solution (see Figure~\ref{fig:ReferenceSolution1D}). This is obtained by applying the a priori matching (Algorithm~\ref{alg:a-priori}) on the uniform grid of the parametric space $\mathcal{M}$ containing 129 points. The matching information is then propagated from left to right along all the 128 subintervals of $\mathcal M$.

\begin{figure}[htpb]
\centering
\includegraphics[width=6cm]{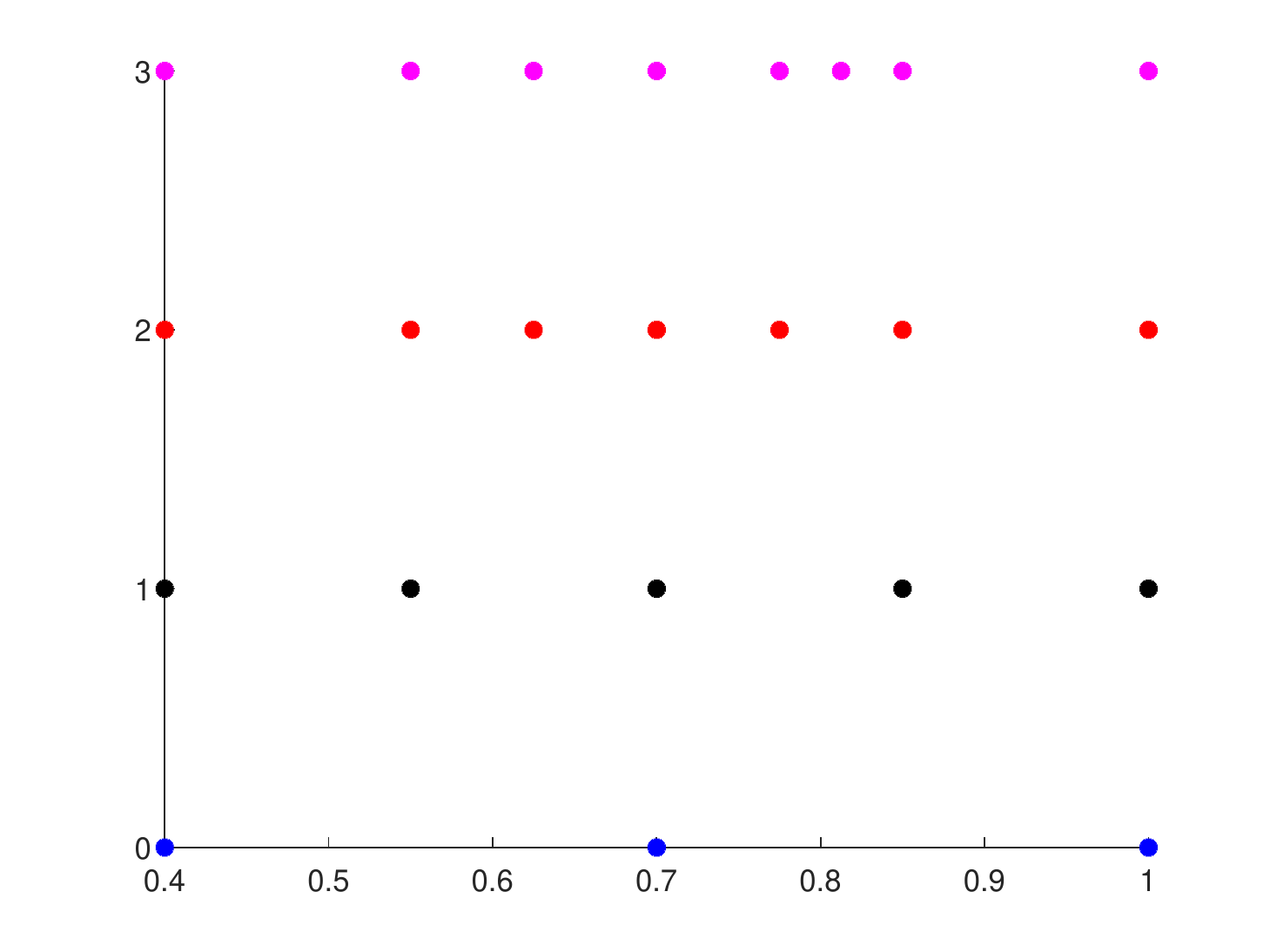}
\caption{Evolution of parametric grid ($x$-axis) as the level increases ($y$-axis).}
\label{fig:GridByLevel}
\end{figure}

\begin{figure}
    \centering
    \subfloat[\centering Level 0 ]{{\includegraphics[width=6cm]{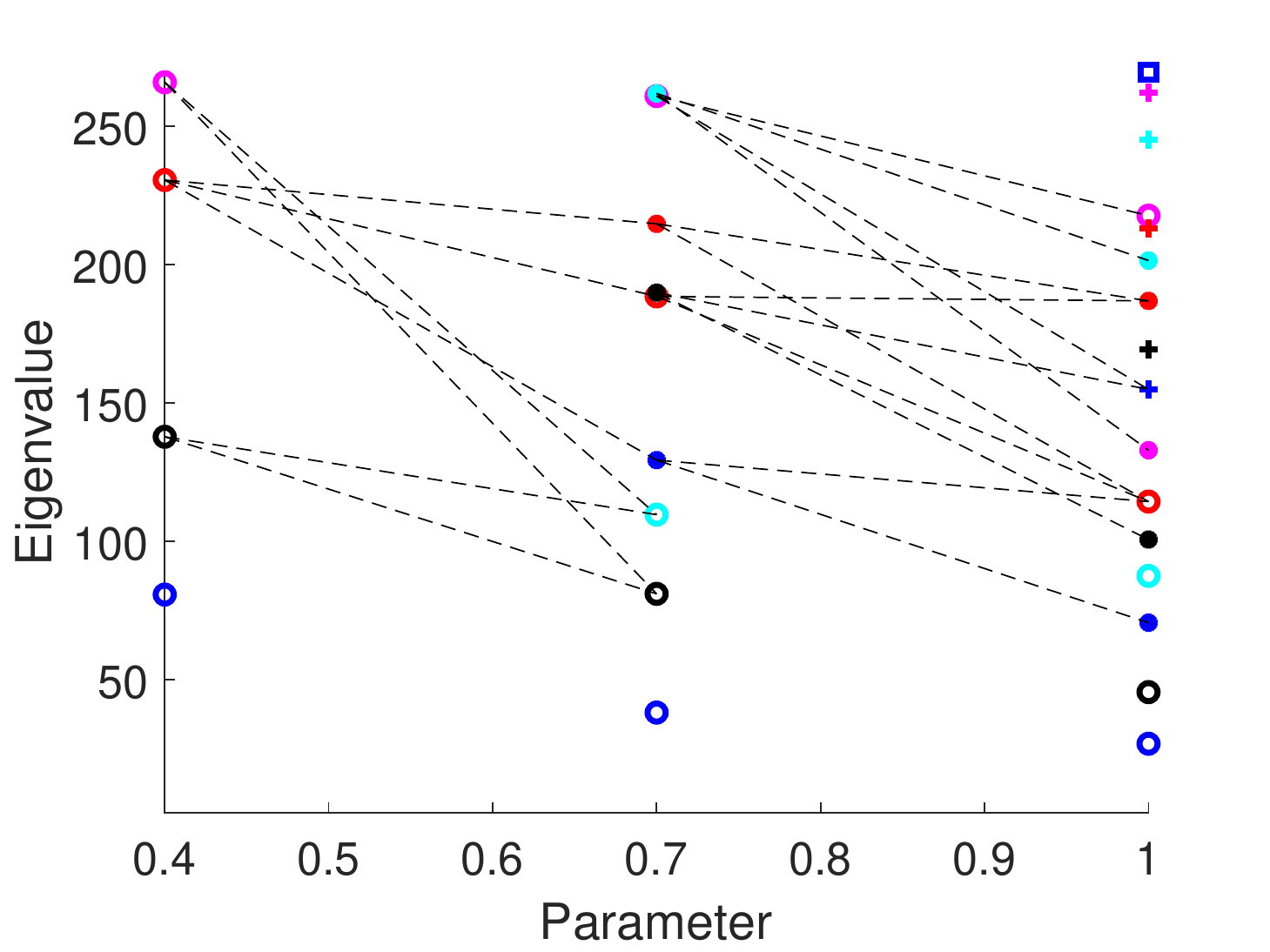}}}
    \qquad
    \subfloat[\centering Level 1 ]{{\includegraphics[width=6cm]{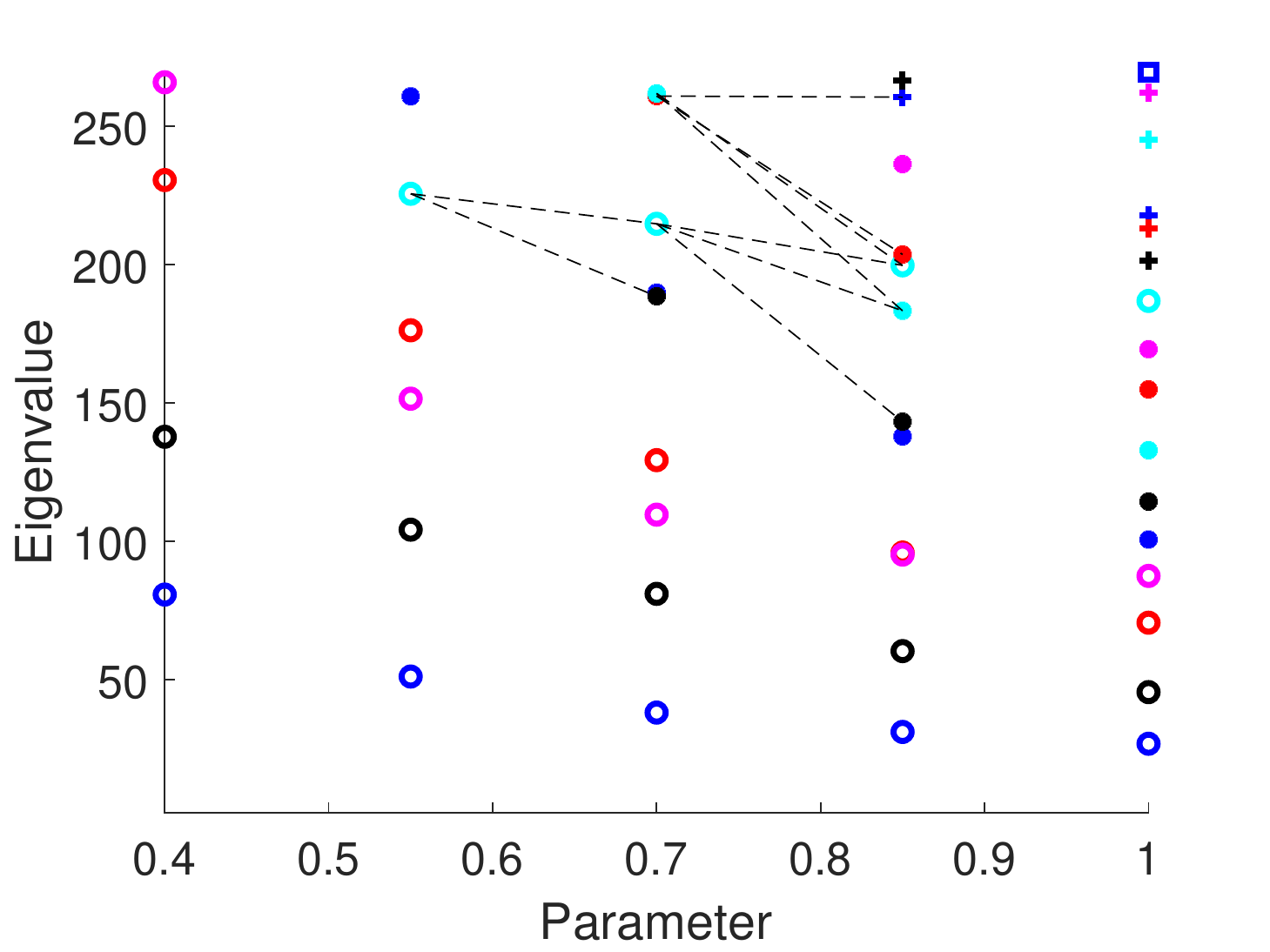}}}
    \qquad
    \subfloat[\centering Level 2 ]{{\includegraphics[width=6cm]{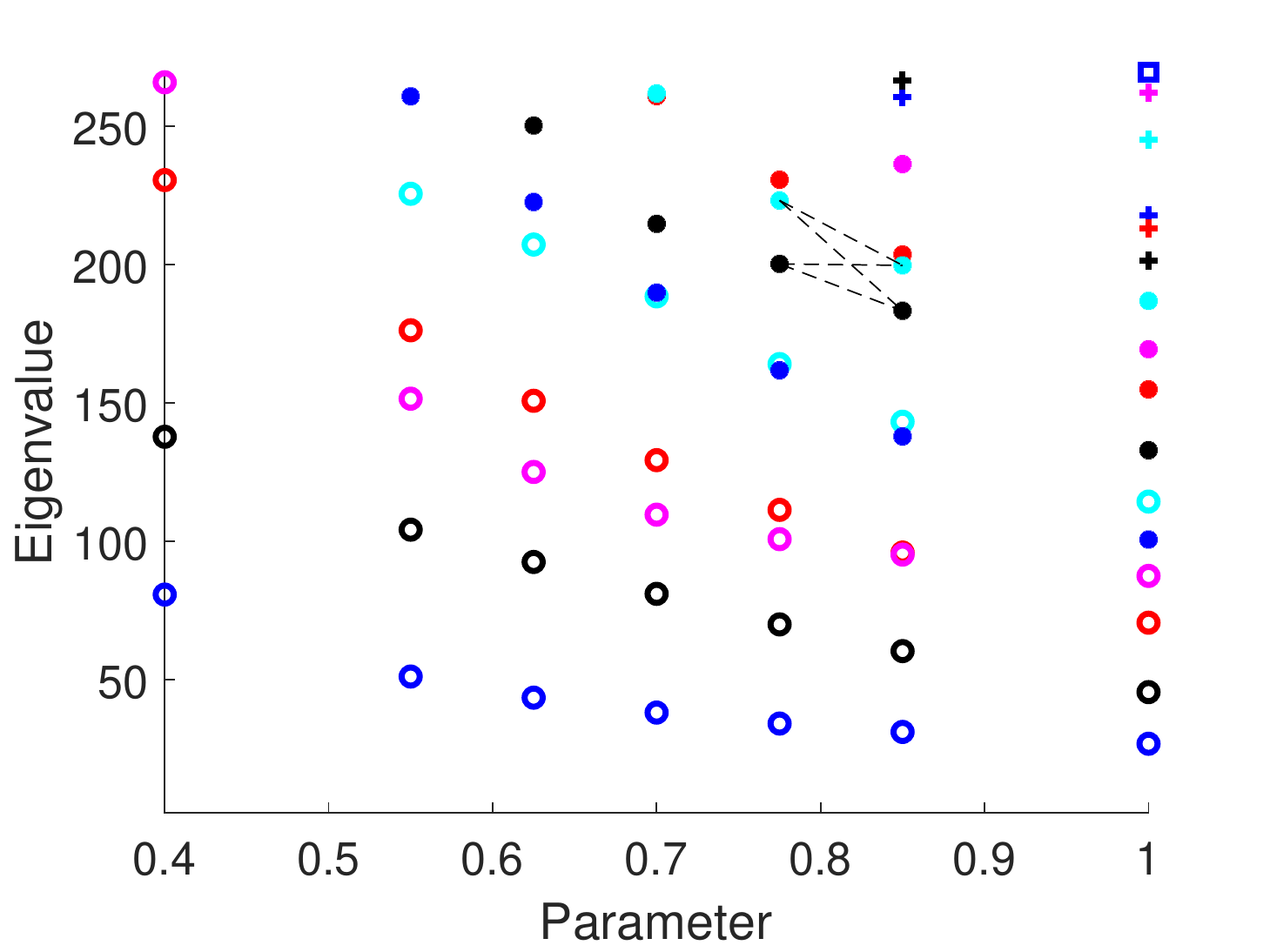}}}
    \qquad
    \subfloat[\centering Level 3 ]{{\includegraphics[width=6cm]{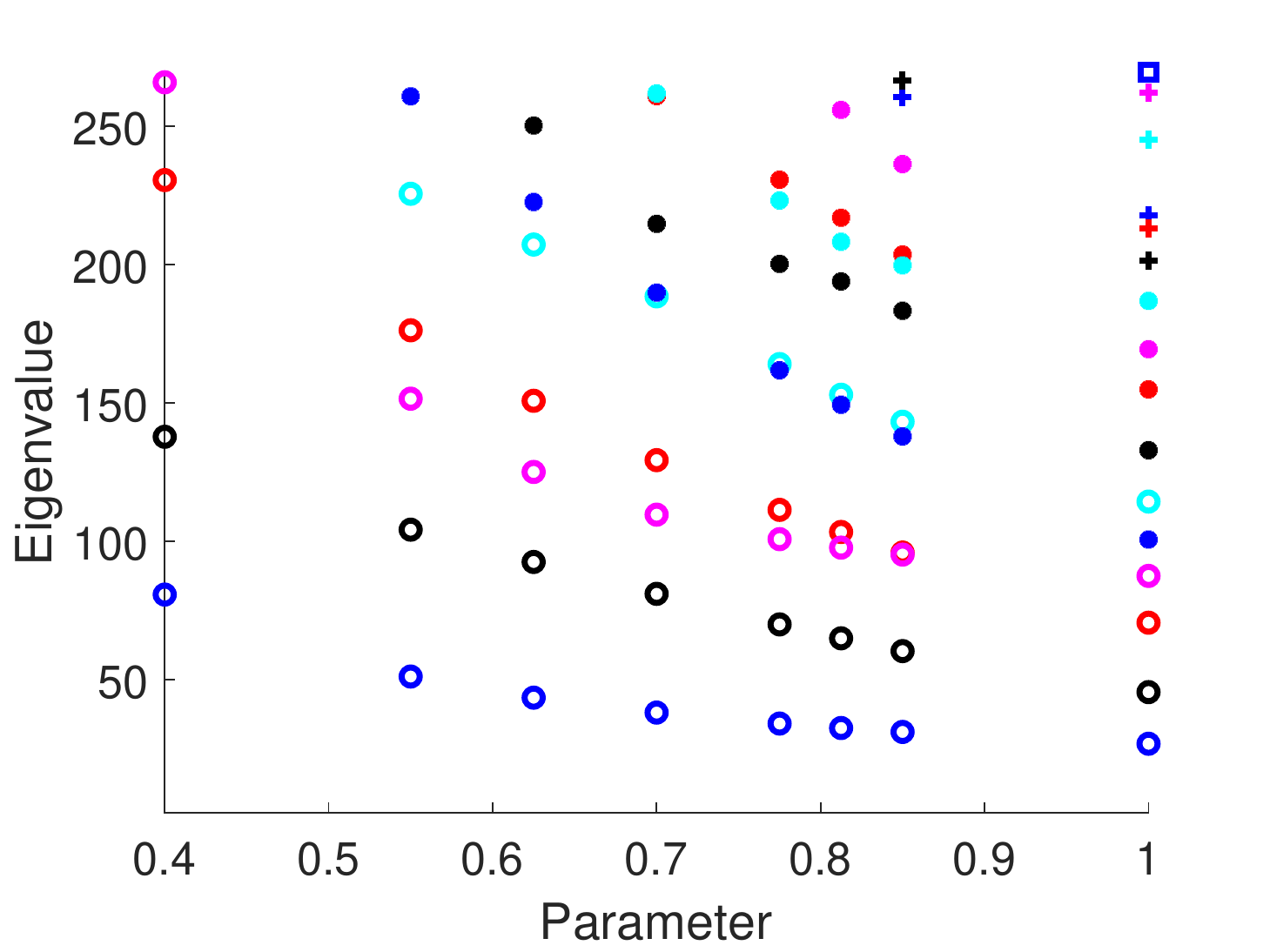}}}
    \caption{Visual summary of the projection matrices.}
    \label{fig:ProjMatrixSummary1D}
\end{figure}

In contrast, we apply the adaptive refinement algorithm presented in Section~\ref{se:Algorithm}. The output is the (coarse) adapted sparse grid of $\mathcal M$ containing enough points to detect all the features of the reference solution. Figure~\ref{fig:GridByLevel} depicts the evolution of the parametric grid as the level increases, whereas in each subfigure of Figure~\ref{fig:ProjMatrixSummary1D} the result of the local checks are represented: matched eigenvalues are plotted using the same color and marker; if the projection matrix suggests that more than one matching is possible, the possible matchings are highlighted by means of black dashed lines.
We detail now the level-by-level procedure, which is summarized in Table~\ref{table:1DError}.

\begin{description}
    \item[Level 0]
    The initial grid is $P^{(0)}=\{0.4,0.7,1\}$ (see the blue points in Figure~\ref{fig:GridByLevel}). In particular, the local checks (both the a priori matching and the a posteriori verification) are performed on the two local subintervals $[0.4, 0.7]$ and $[0.7, 1]$, which are both marked for further refinement.
    
    \item[Level 1]
    The grid at level one contains five points $P^{(1)}=\{0.4,0.55,0.7,0.85,1\}$ (see the black dots in Figure~\ref{fig:GridByLevel}) and it is given by $P^{(1)}=P^{(0)}\cup P^{(1)}_\delta$, with $P^{(1)}_\delta=\{0.55,0.85\}$. Following Algorithm~\ref{alg:refinement}, the local checks are performed on the four subintervals $[0.4,0.55]$, $[0.85,1]$, $[0.55,0.7]$ and $[0.7, 0.85]$. Only the last two subintervals are marked for further refinement.

    \item[Level 2]
    The grid at level 2 is $P^{(2)}=\{0.4,0.55,0.625,0.7,0.775,0.85,1\}$ (see the red dots in Figure~\ref{fig:GridByLevel}) and it is given by $P^{(2)}=P^{(1)}\cup P^{(2)}_\delta$, with $P^{(2)}_\delta=\{0.625,0.775\}$. Local checks are then performed on the following subintervals: $[0.55,0.625]$, $[0.625,0.7]$, $[0.7,0.775]$, and $[0.775, 0.85]$. Only one subinterval is marked for further refinement, namely $[0.775,0.85]$. 
    
    \item[Level 3]
    The grid at level 3 is $P^{(3)}=\{0.4,0.55,0.625,0.7,0.775,0.8125,0.85,1\}$ (see the magenta dots in Figure~\ref{fig:GridByLevel}) and it is given by $P^{(3)}=P^{(2)}\cup P^{(3)}_\delta$, with $P^{(3)}_\delta=\{0.8125\}$. The local checks on the local subintervals $[0.775,0.8125]$ and $[0.8125,.085]$ are performed, and none of the two is marked for refinement. As a consequence, the adaptive algorithm terminates.
\end{description}

\begin{figure}
    \centering
    \subfloat[\centering $\lambda = 200.3$, $\mu = 0.775$ ]{{\includegraphics[width=6cm]{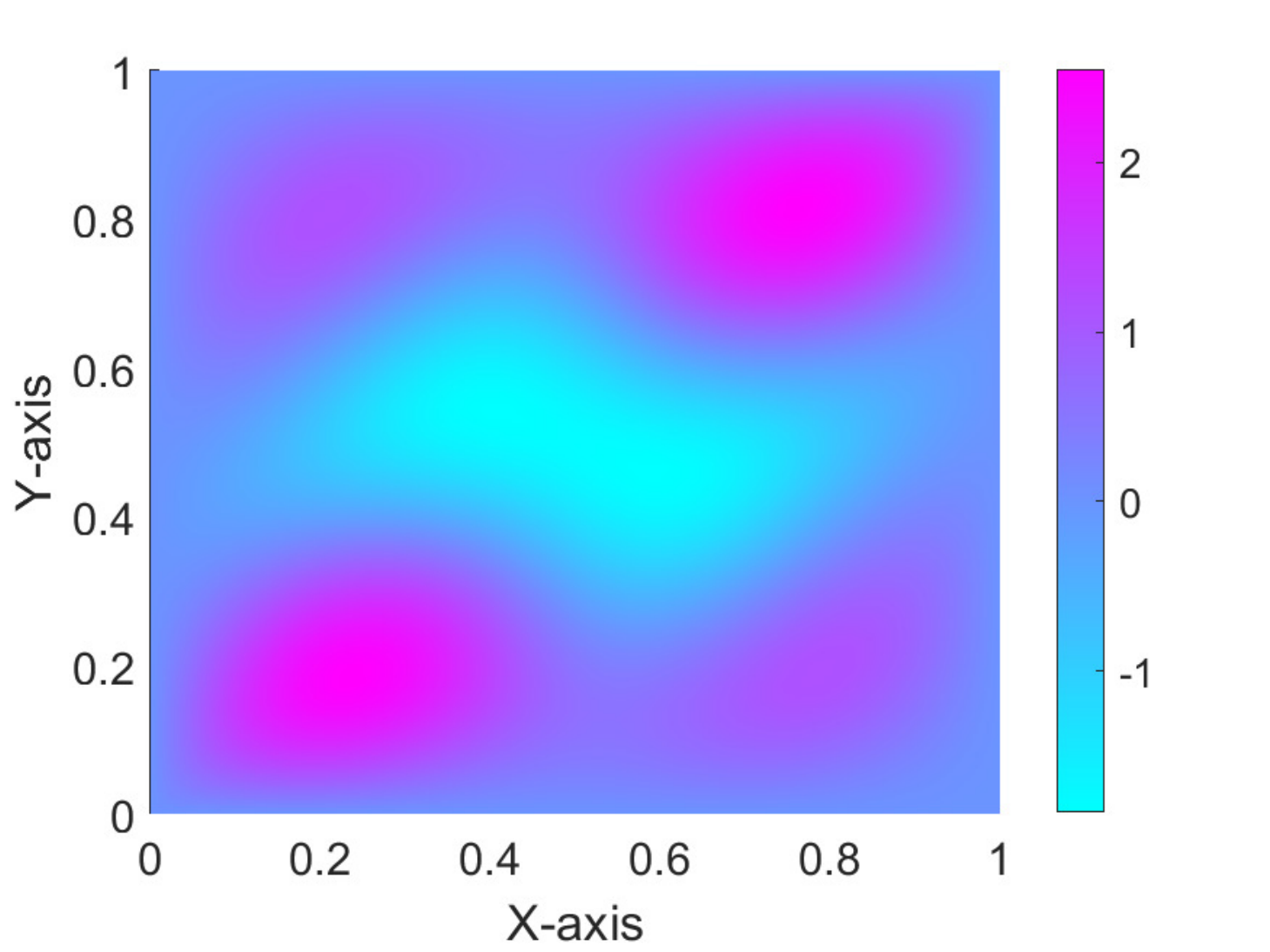}}}
    \qquad
    \subfloat[\centering $\lambda = 223.23$, $\mu = 0.775$ ]{{\includegraphics[width=6cm]{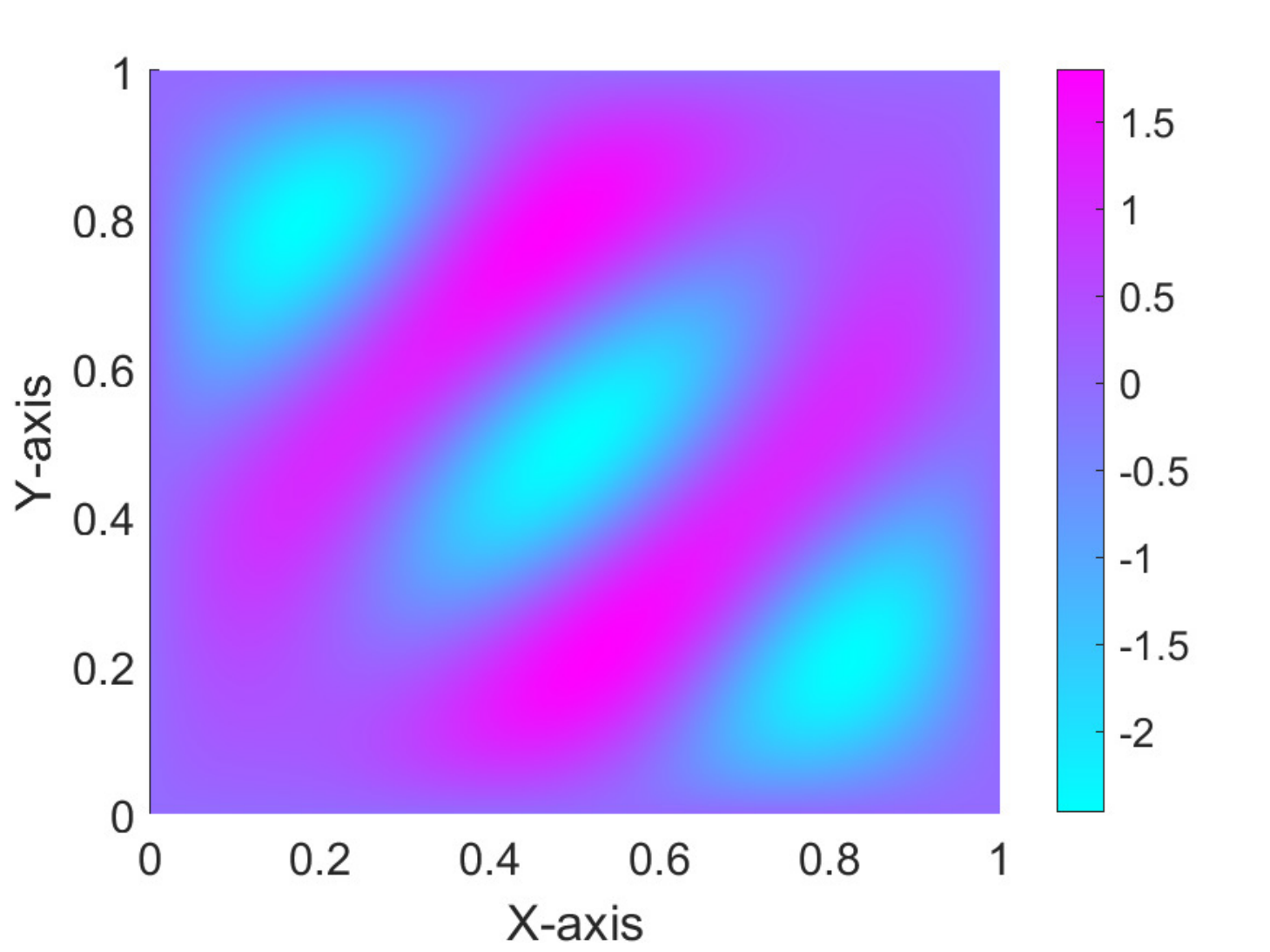}}}
    \qquad
    \subfloat[\centering $\lambda = 183.369$, $\mu = 0.85$ ]{{\includegraphics[width=6cm]{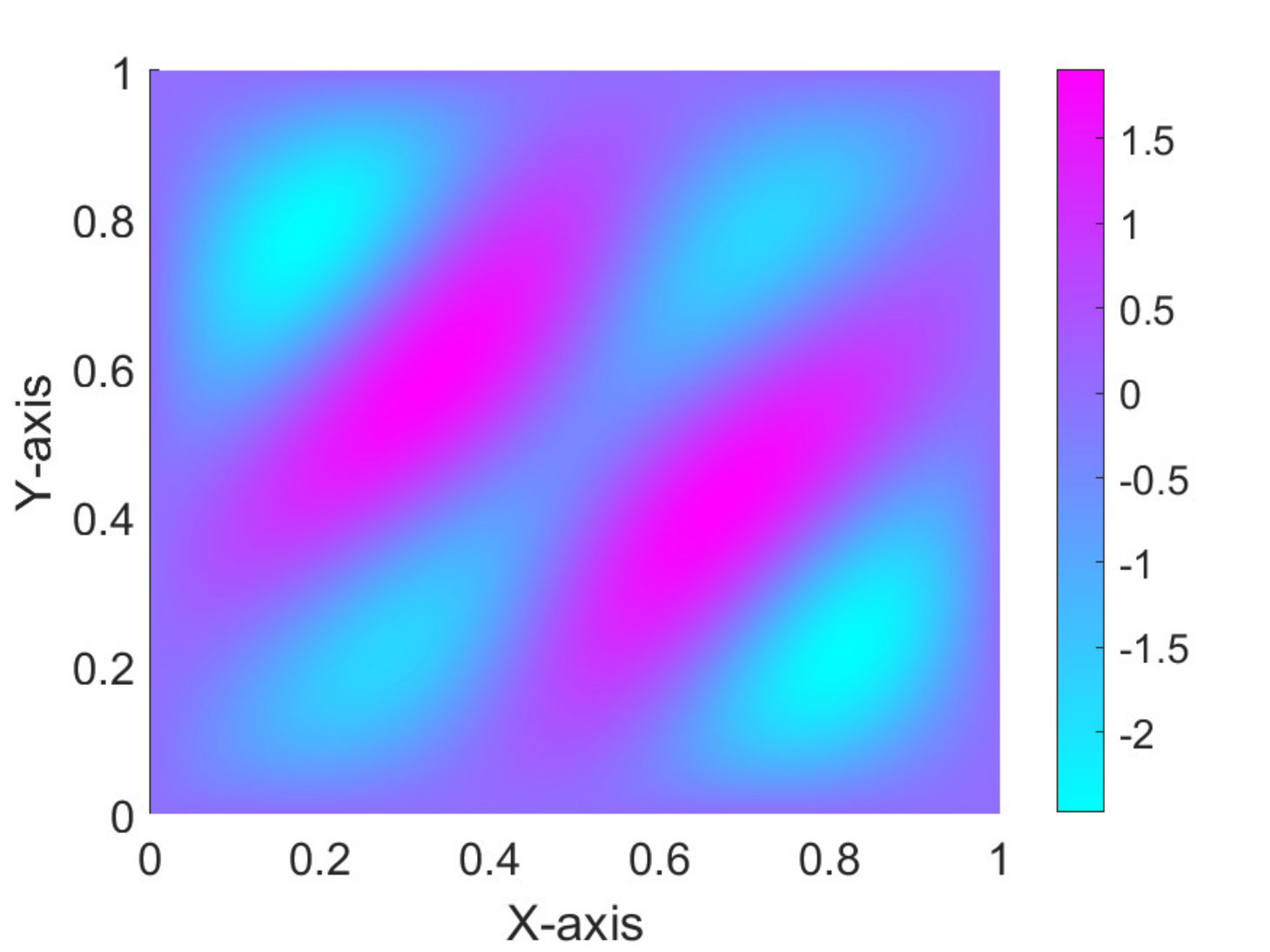}}}
    \qquad
    \subfloat[\centering $\lambda = 199.786$, $\mu = 0.85$ ]{{\includegraphics[width=6cm]{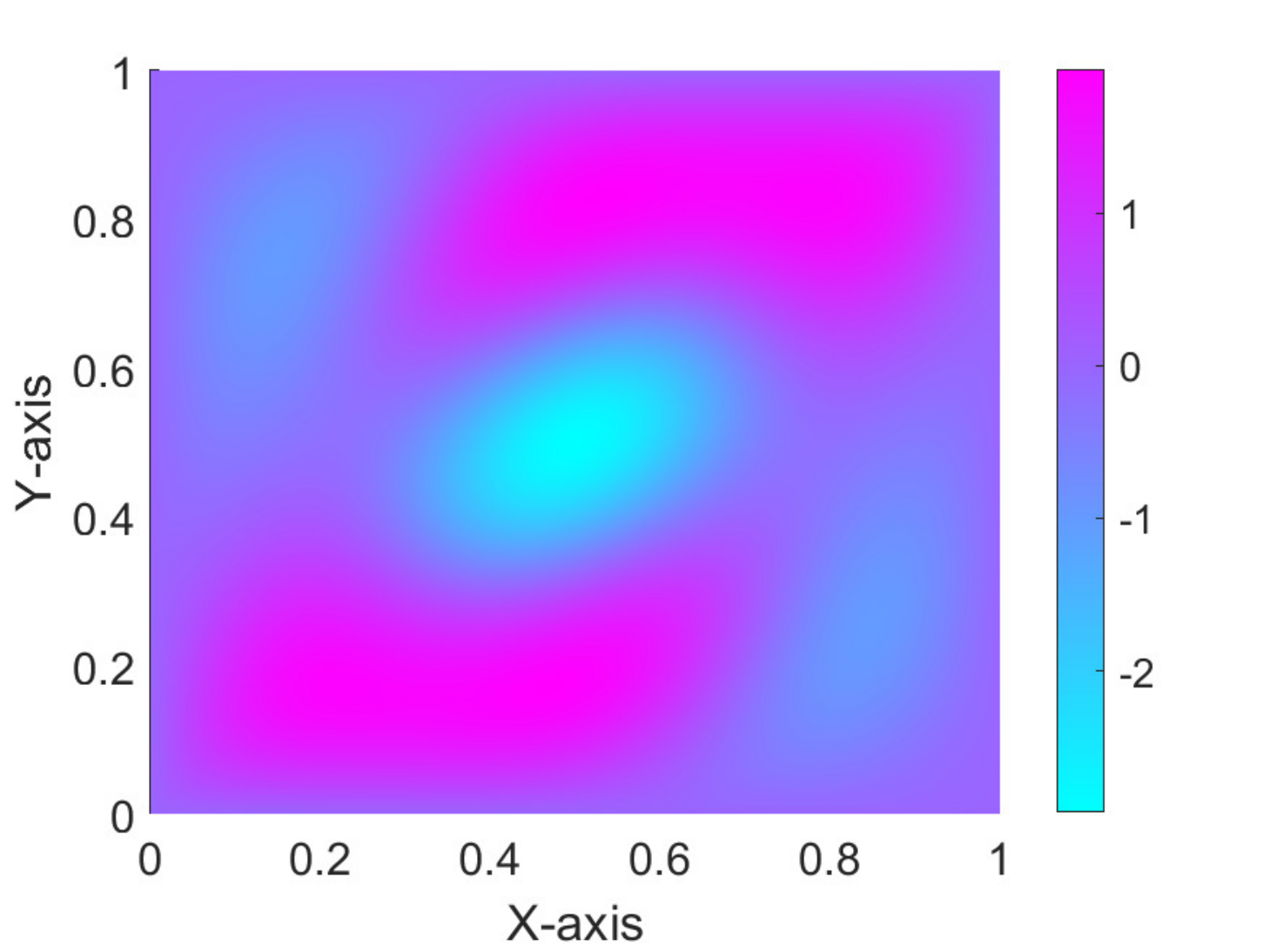}}}
    \caption{Comparison between four non-orthogonal eigenfunctions.}
    \label{fig:EigenfunctionComparison}

\end{figure}

It is worth mentioning that, even though the subinterval $[0.775,0.85]$ was marked for refinement at level 2, the a priori matching was correct. To explain the extra refinement, we examine the projection matrix associated with this subinterval, more precisely the $2\times 2$ sub-matrix with entries corresponding to the black and cyan eigenvalues:
$$
\begin{bmatrix}
\Pi_{7,7}  &  \Pi_{7,8} \\
\Pi_{8,7} &  \Pi_{8,8}
\end{bmatrix} 
=
\begin{bmatrix}
0.721  &  0.682 \\
0.680 &  0.714
\end{bmatrix}.
$$
The off-diagonal terms $\Pi_{7,8}$ and $\Pi_{8,7}$ are non-zero due to a non-negligible overlap between the corresponding eigenfunctions, which are depicted in Figure~\ref{fig:EigenfunctionComparison}. Therefore, we can identify two causes for the a posteriori verification to suggest the refinement: (i) the ordering of the eigenvalues proposed by the a priori matching is incorrect; (ii) the eigenfunctions are not orthogonal. The second situation typically occurs in the presence of small gaps between the eigenvalue hypersurfaces. Recall that lines 20-25 in Algorithm~\ref{alg:a-posteriori2} are devoted to case (ii). In this example, no clusters are identified, and the algorithm consequently terminates.

\begin{figure}
    \centering
    \subfloat[\centering Level 0 ]{{\includegraphics[width=6cm]{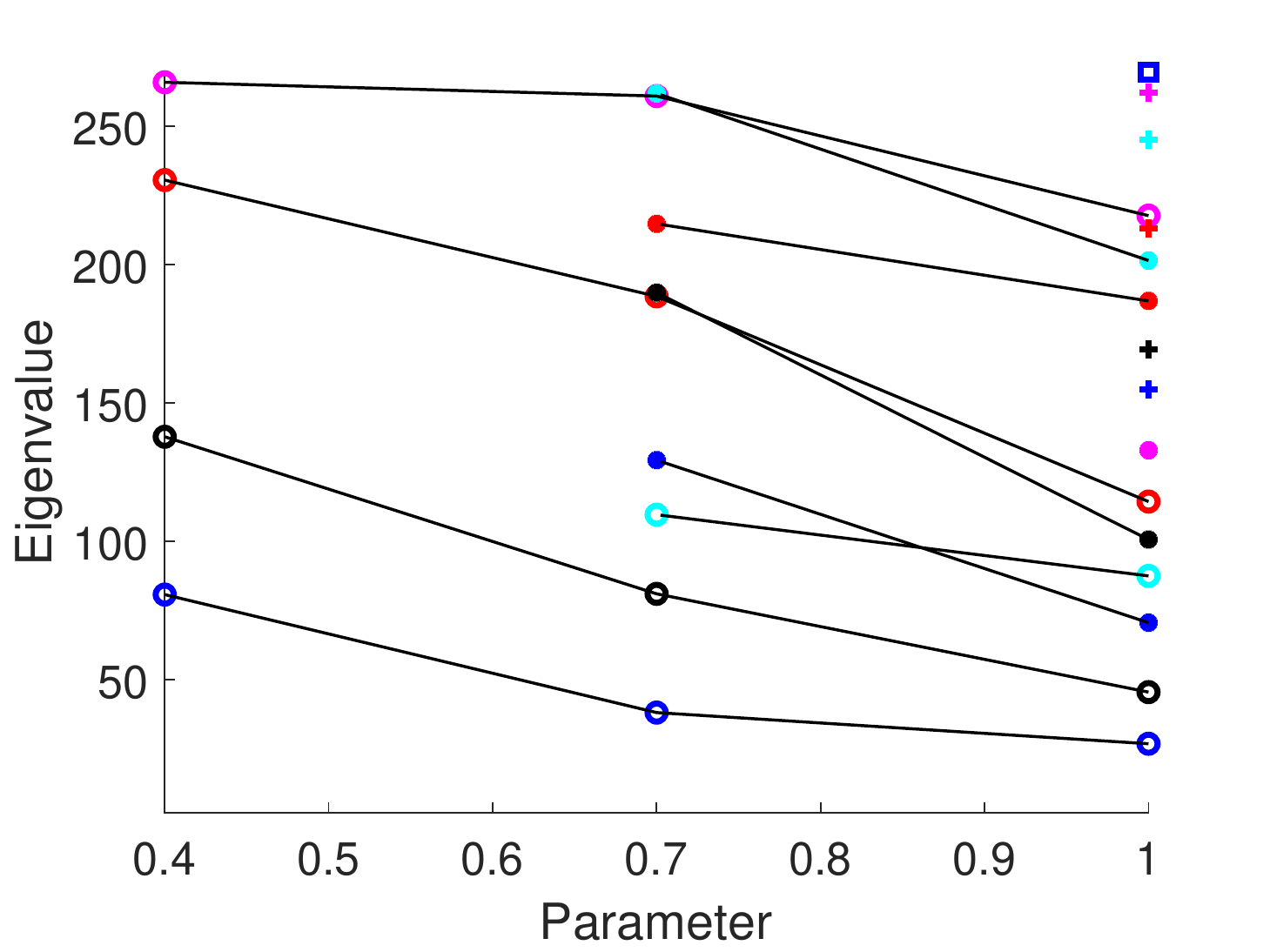}}}
    \qquad
    \subfloat[\centering Level 1 ]{{\includegraphics[width=6cm]{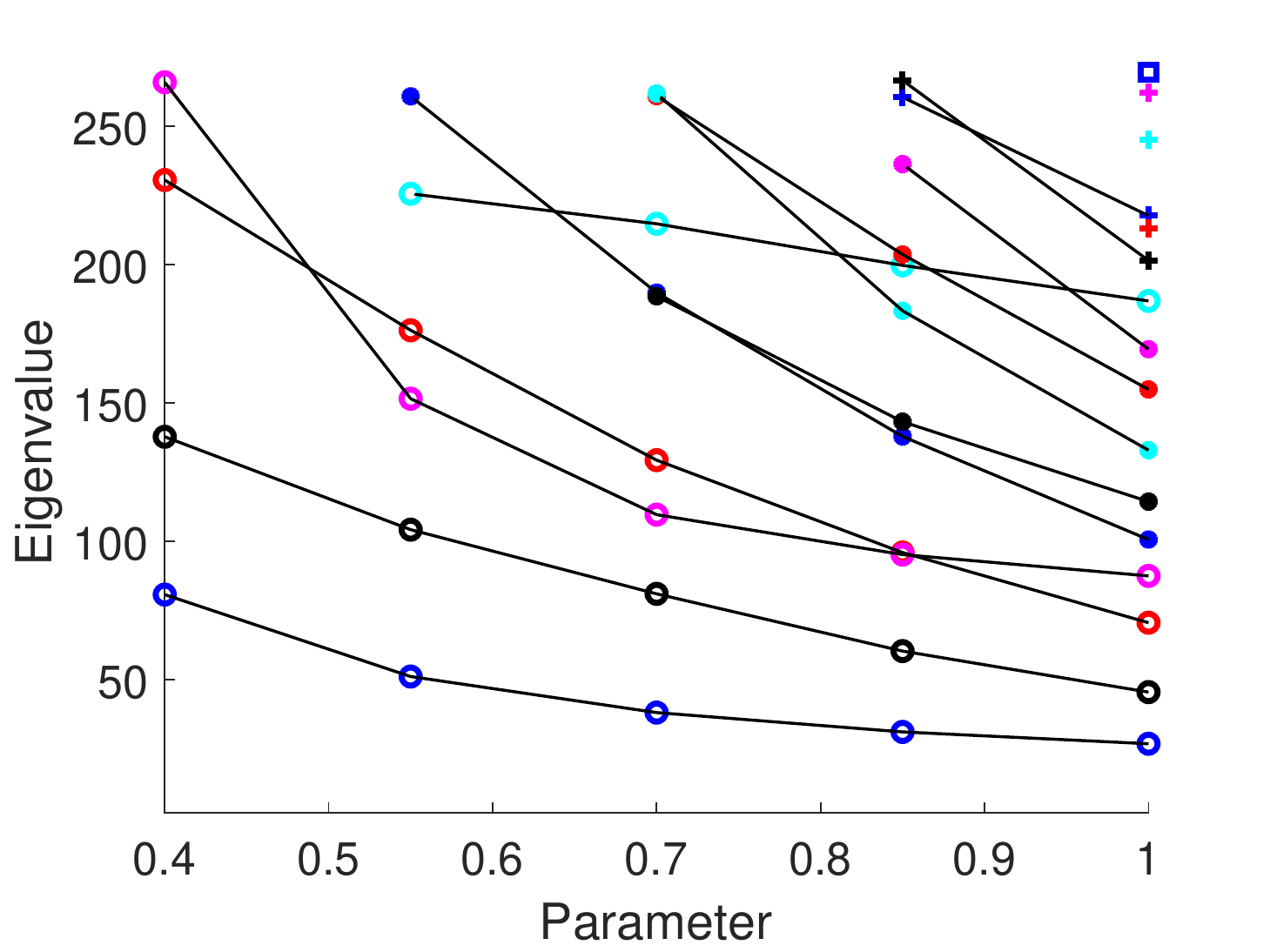}}}
    \qquad
    \subfloat[\centering Level 2 ]{{\includegraphics[width=6cm]{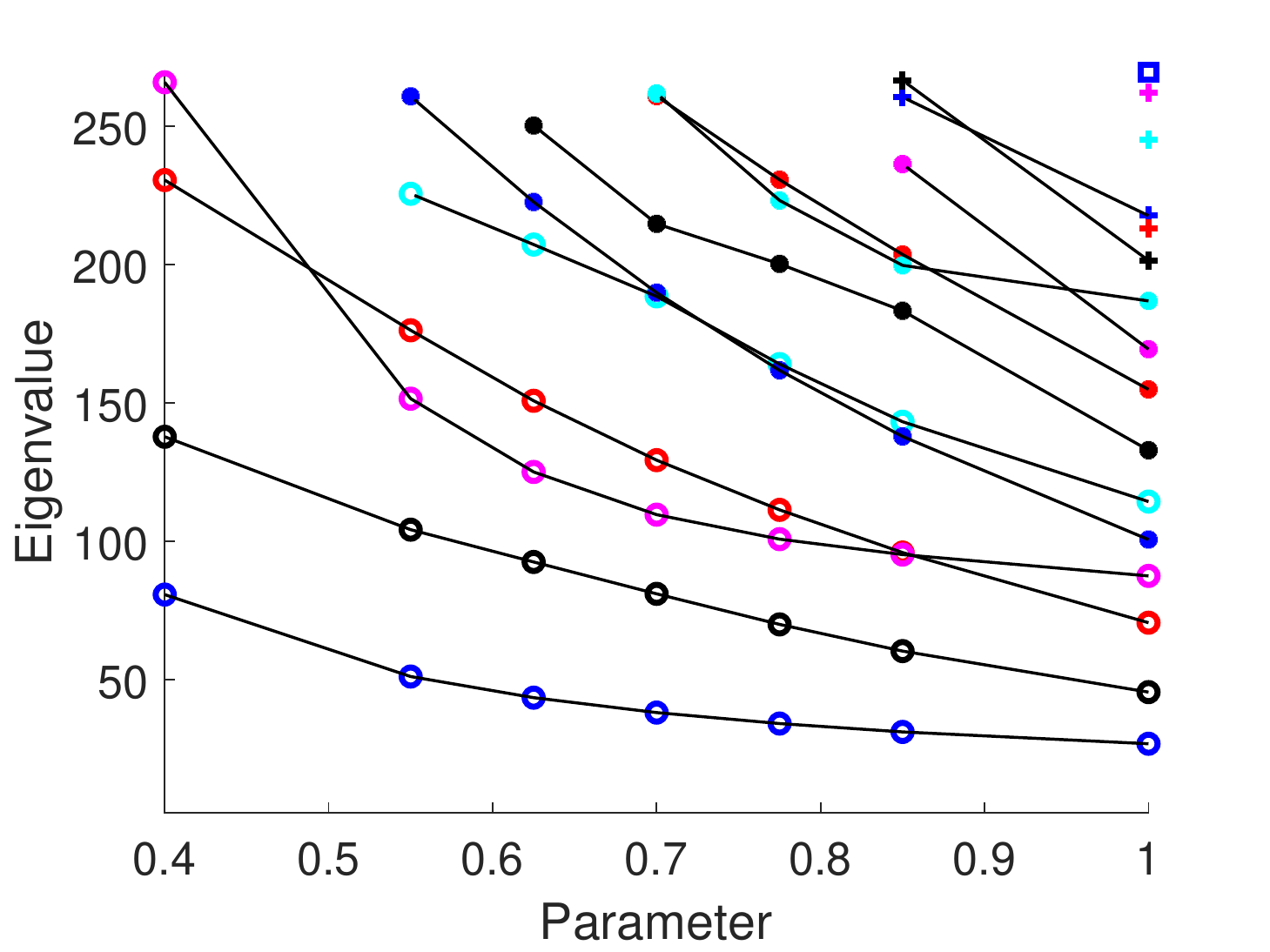}}}
    \qquad
    \subfloat[\centering Level 3 ]{{\includegraphics[width=6cm]{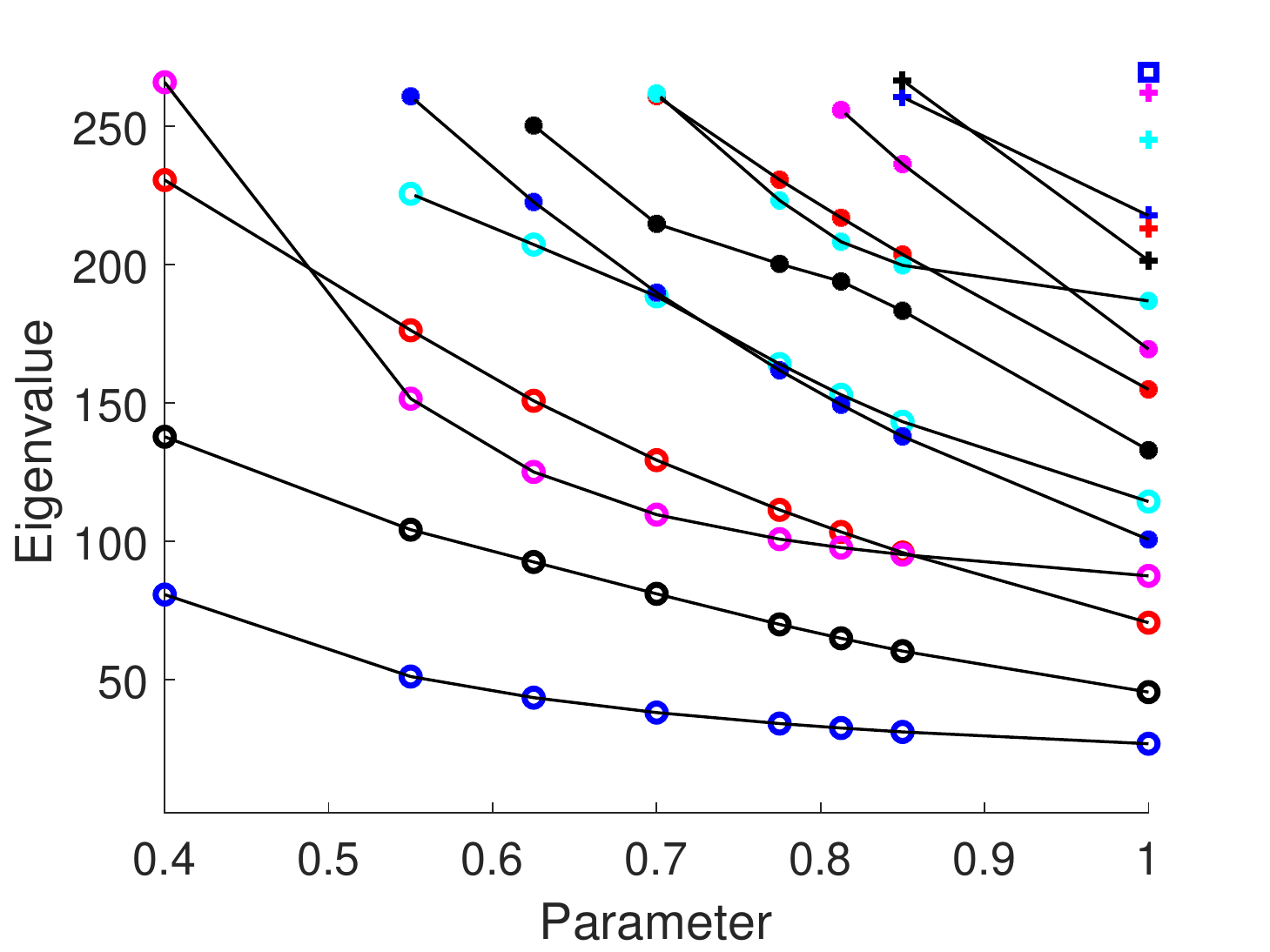}}}
    \caption{Output of each level for the 1D example.}
    \label{fig:EigenvaluesByLevel}

\end{figure}

Once the adaptive algorithm has terminated, the collected information can be exploited to construct a surrogate for the map $\mu\mapsto\lambda(\mu)$. Between the various possibilities, we propose to construction this surrogate simply by piecewise linear interpolation. The result is depicted in Figure~\ref{fig:EigenvaluesByLevel}. 

\begin{table}
\begin{center}
\begin{tabular}{ |p{1cm}|p{2cm}|p{2.5cm}|p{2cm}|p{2.5cm}|  }
\hline
\multicolumn{5}{|c|}{Error By Level} \\
\hline
Level & Total no. of points & No. of wrongly matched points & No. of subintervals & No. of uncertified subintervals \\
\hline
0 & 3 & 2 & 2 & 2 \\
\hline
1  & 5 & 3 & 4 & 2 \\
\hline
2  & 7 & 0 & 4 & 1 \\
\hline
3  & 8 & 0 & 2 & 0 \\
\hline
\end{tabular}
\end{center}
\caption{Error table for the 1D example.}
\label{table:1DError}
\end{table}

%%%%%%%%%%%%%%%%%%%%%%%%%%%%%%%%%%%%%%%%%%%%%%%%%%%%%%%%%%%%%%%%%%%%%%%%%%%%%%%%%%
\subsection{A 2D example}
\label{se:2D}
%%%%%%%%%%%%%%%%%%%%%%%%%%%%%%%%%%%%%%%%%%%%%%%%%%%%%%%%%%%%%%%%%%%%%%%%%%%%%%%%%%

In this numerical test we consider the two dimensional parameter space $\mathcal M=[0.8, 1.05] \times [0.8, 1.05]$, and we take the (positive definite) diffusion matrix 
$$
c(\mu) = c(\mu_1,\mu_2) = 
\begin{pmatrix}
\mu_1^{-2} & 0.8 \mu_2^{-1} \\
0.8 \mu_2^{-1} & \mu_2^{-2}
\end{pmatrix}.
$$

\begin{figure}[htpb]
\centering
\includegraphics[width=6cm]{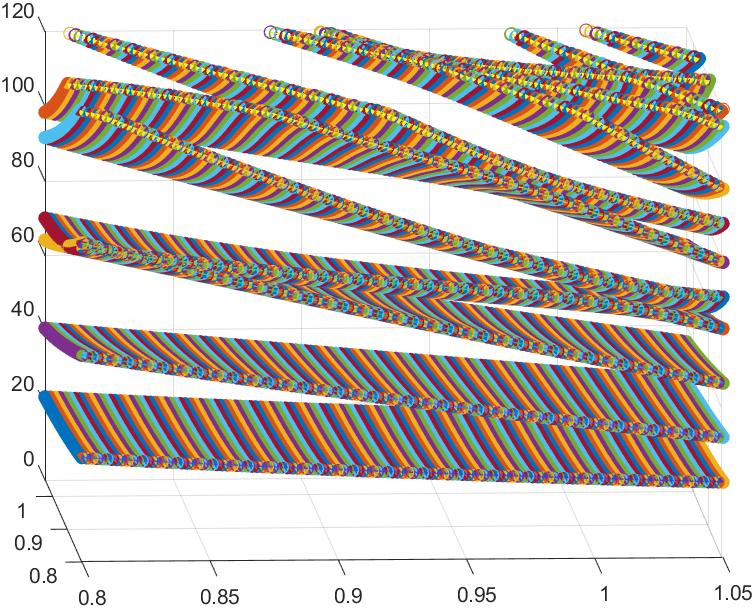}
\includegraphics[width=6cm]{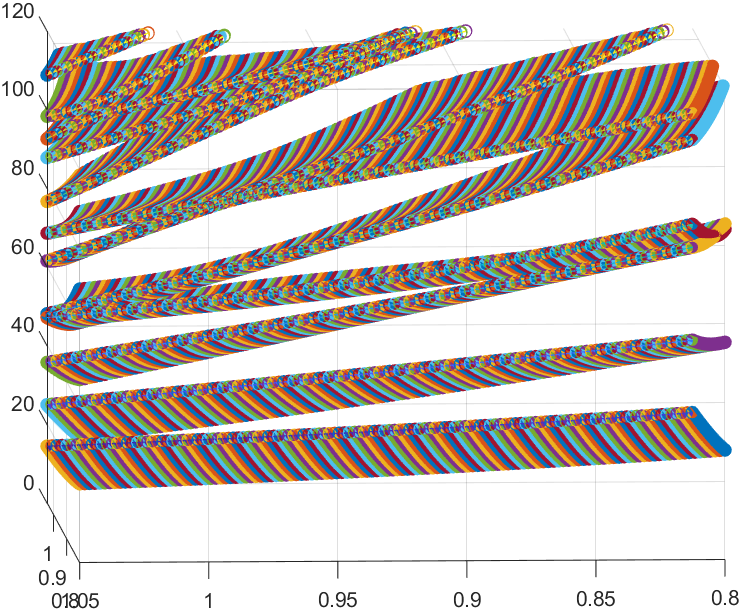}
\caption{Reference solution for the 2D example; All hypersurfaces.}
\label{fig:2D_reference_sol}
\end{figure}

The computation of the reference solution relies on the uniform tensor product grid of $\mathcal M$ containing $129$ points (see Figure~\ref{fig:2D_reference_sol} for the eigenvalue hypersurfaces). For eigenvalue problems depending on two parameters, the computation of the reference solution is still affordable. However, for increasing dimension $d$ of the parametric space, such an approach entails prohibitive computational costs, becoming out of reach and extremely more expensive than the proposed sparse grid-based approach.
In the following discussion, we focus on the two important features of the reference solution, namely, crossings and small gaps of the eigenvalue hypersurfaces.

\begin{figure}[htpb]
    \centering
    \subfloat[\centering Level 0 ]{{\includegraphics[width=6cm]{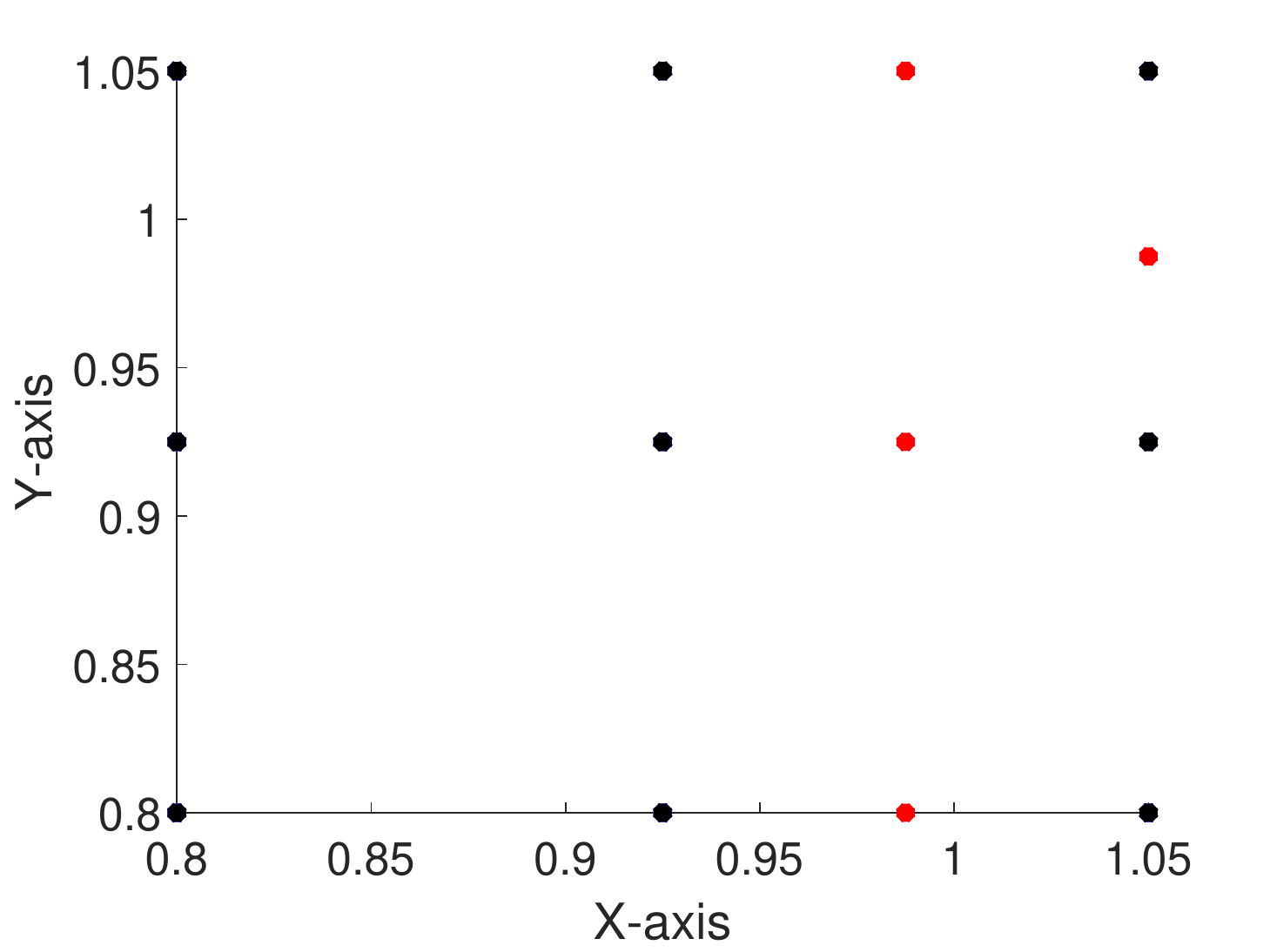}}}
    \qquad
    \subfloat[\centering Level 1 ]{{\includegraphics[width=6cm]{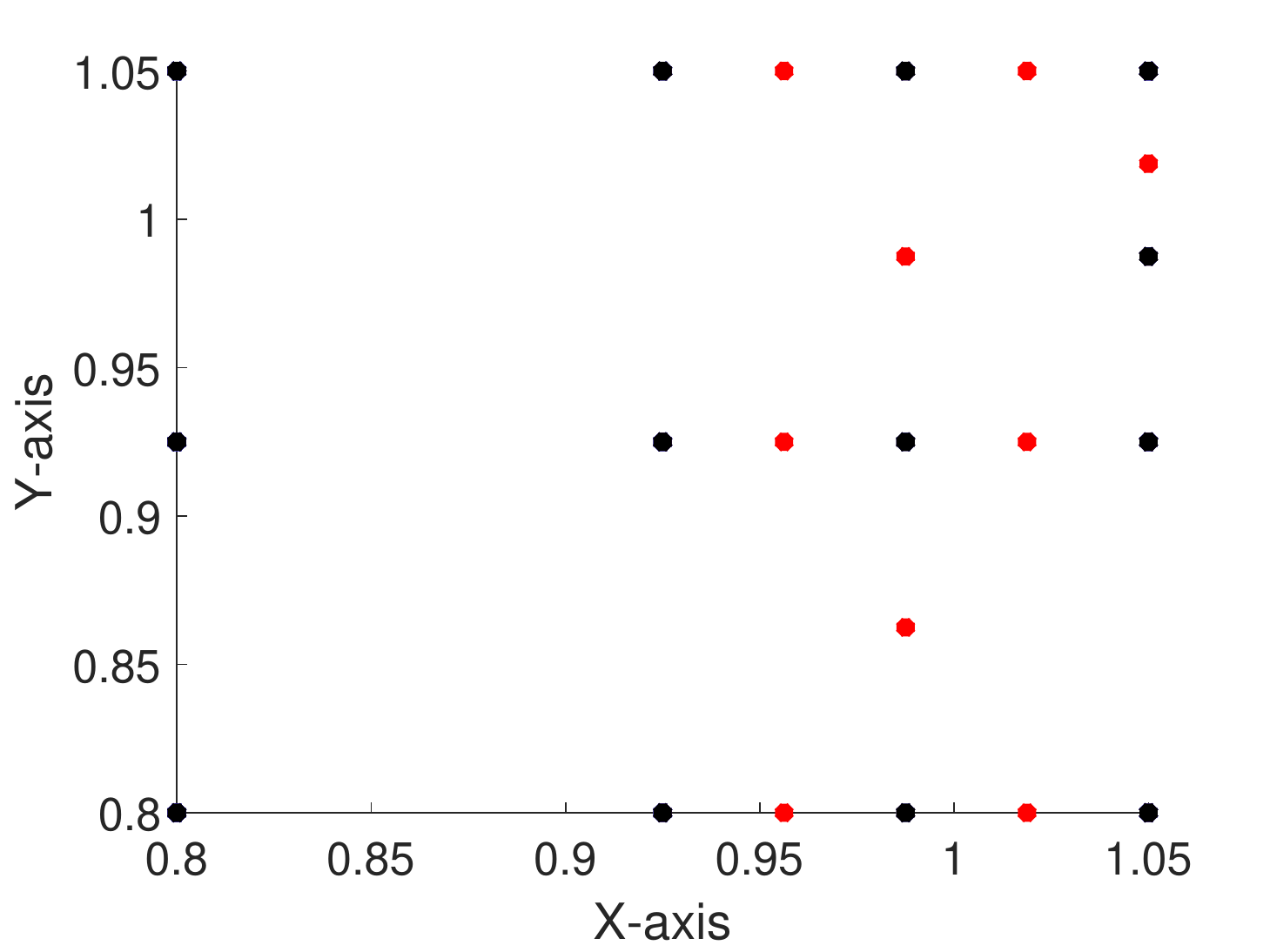}}}
    \qquad
    \subfloat[\centering Level 2 ]{{\includegraphics[width=6cm]{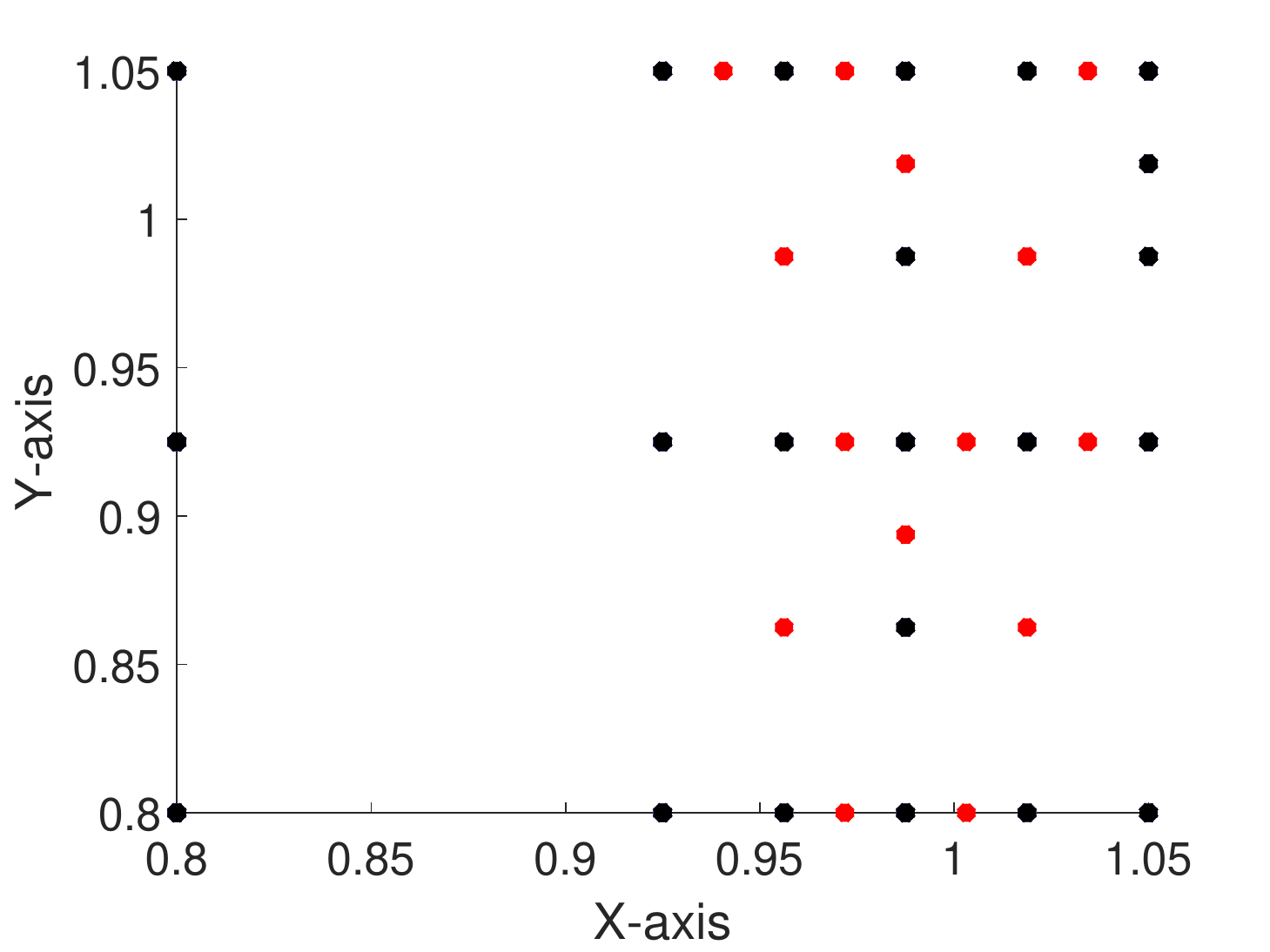}}}
    \qquad
    \subfloat[\centering Level 3 ]{{\includegraphics[width=6cm]{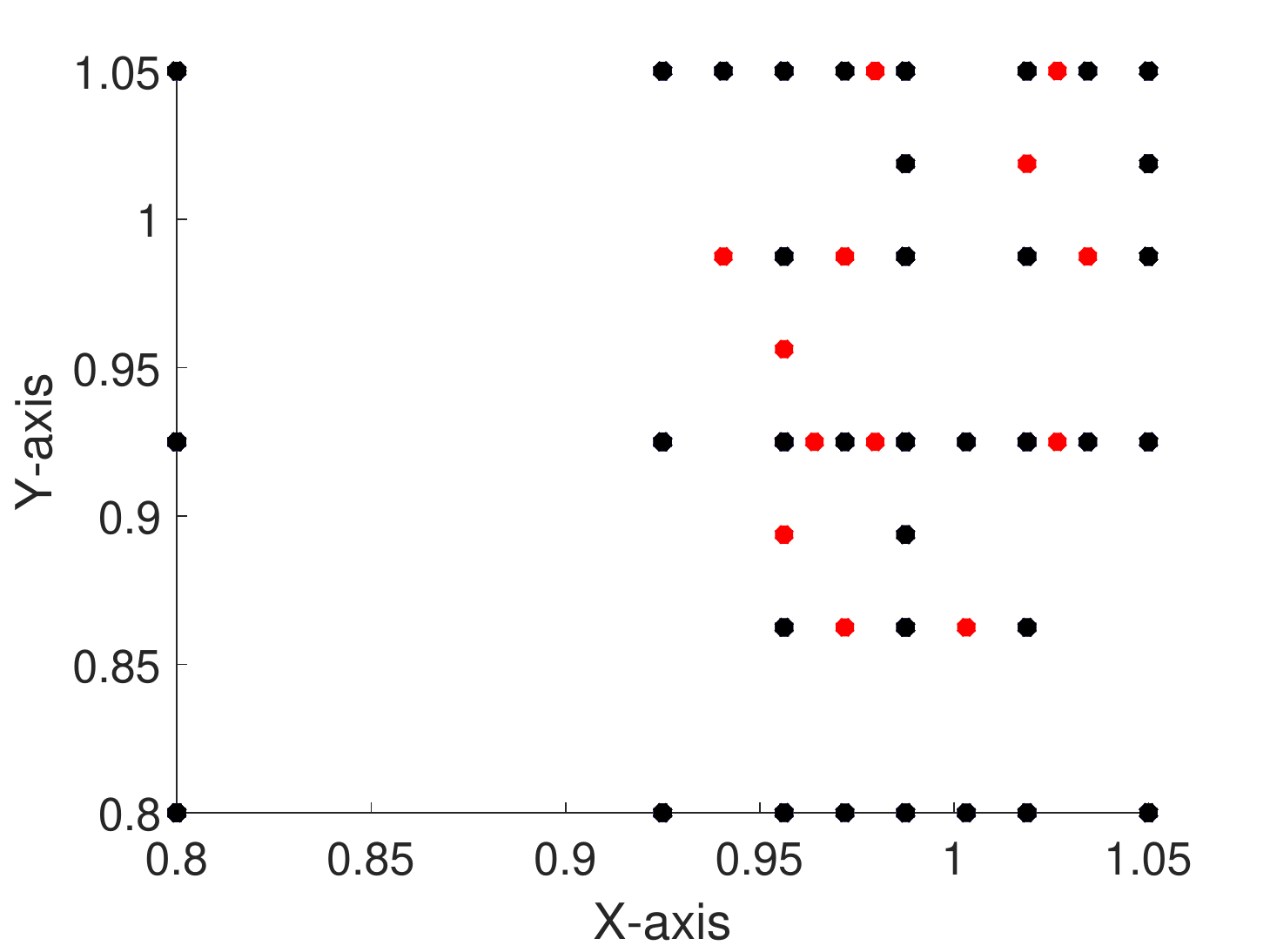}}}
    \qquad
    \subfloat[\centering Level 4 ]{{\includegraphics[width=6cm]{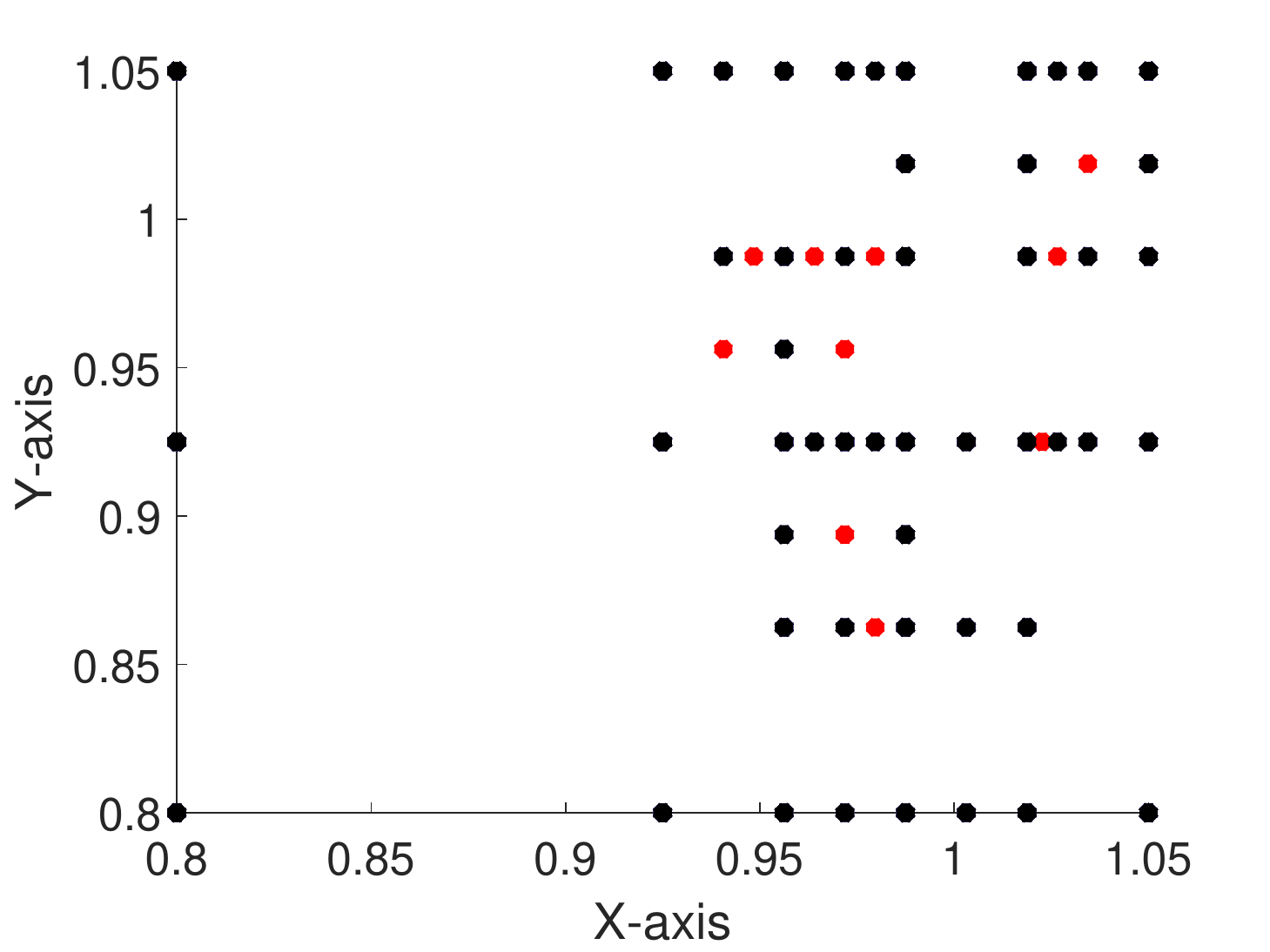}}}
    \qquad
    \subfloat[\centering Level 5 ]{{\includegraphics[width=6cm]{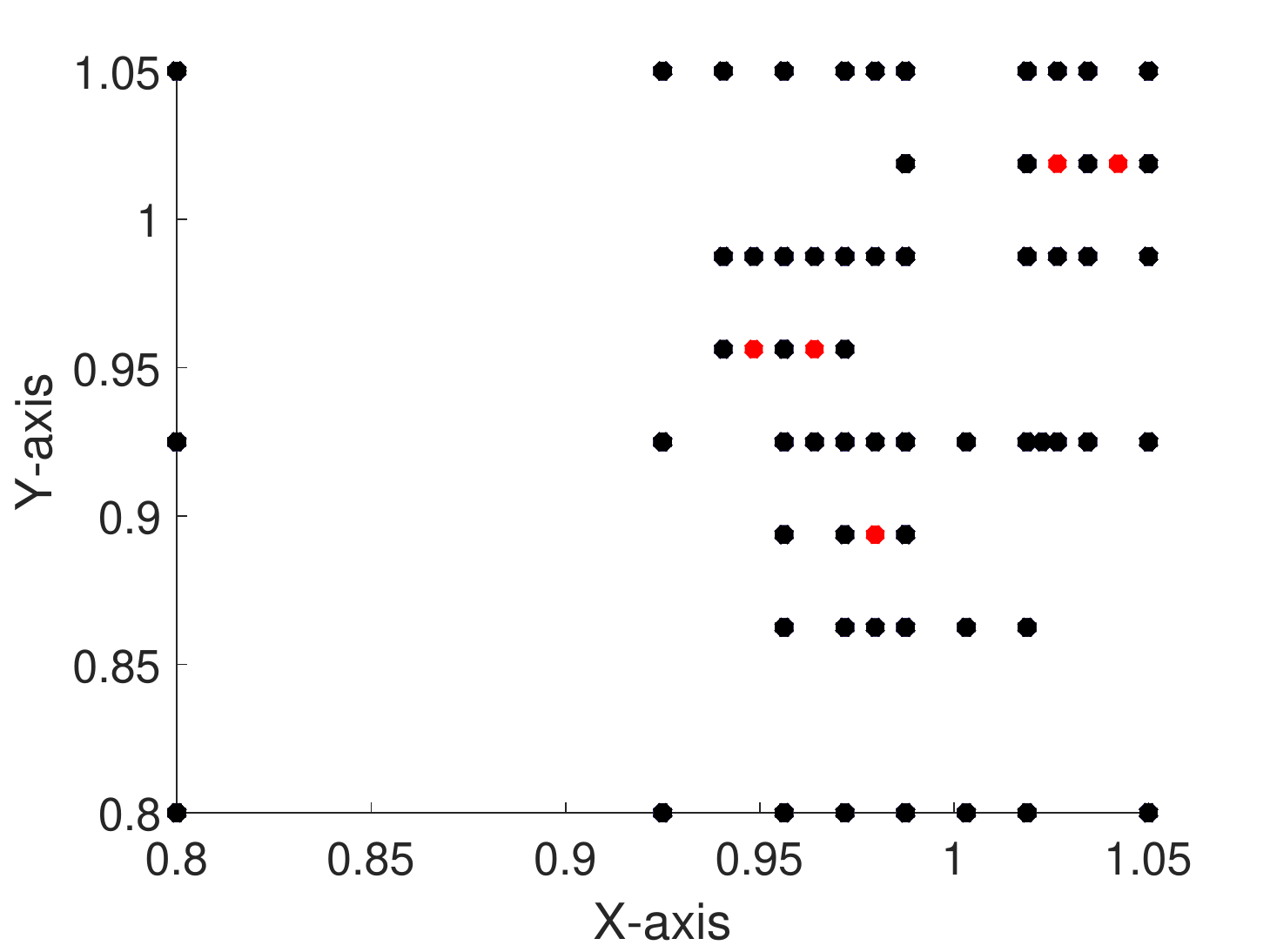}}}
    \qquad
    \subfloat[\centering Level 6 ]{{\includegraphics[width=6cm]{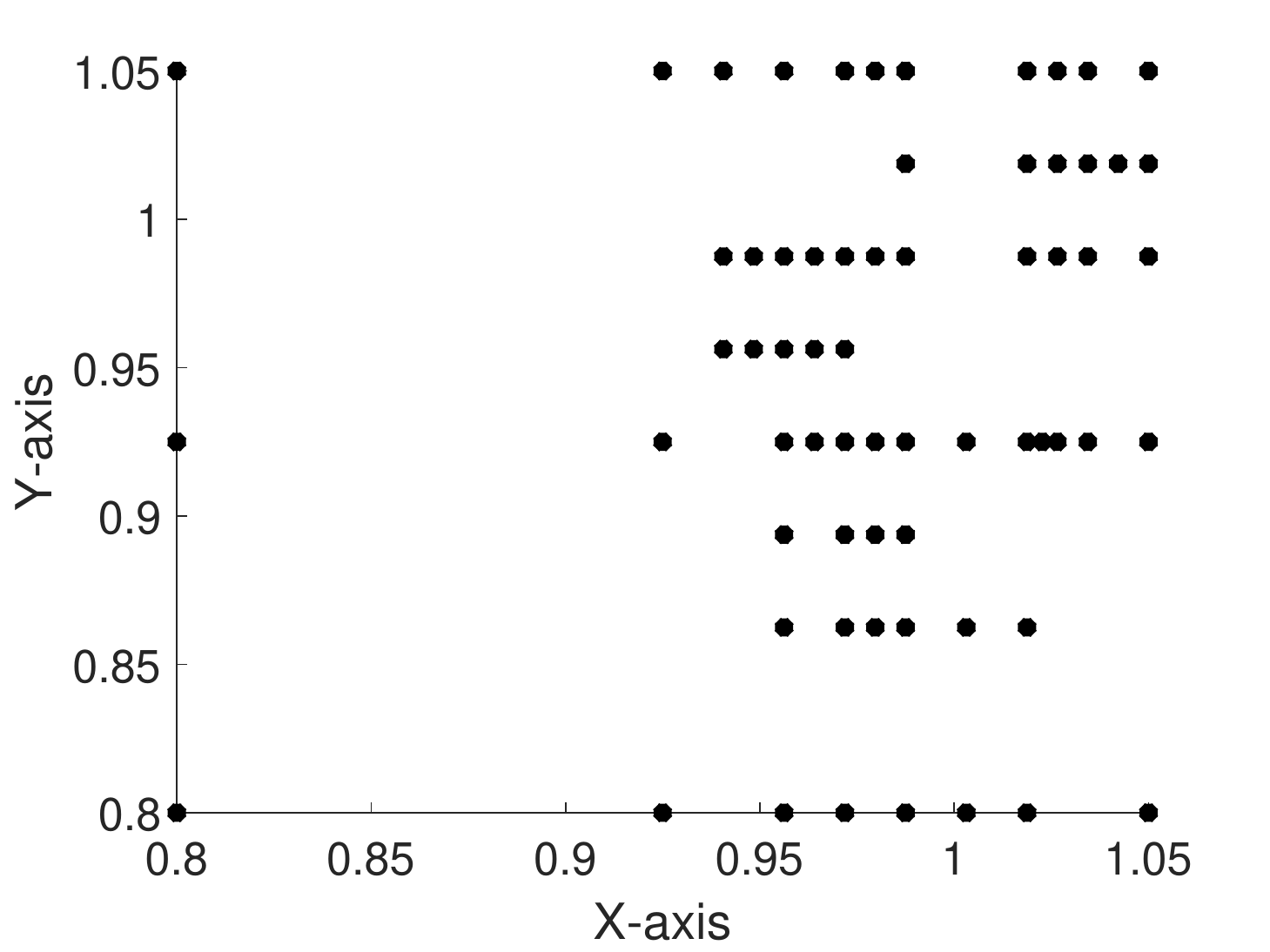}}}
    \caption{Evolution of grid in the 2D example.}
    \label{fig:2D_grid_level_by_level}
\end{figure}

\begin{figure}
    \centering
    \subfloat[\centering Hypersurfaces 4,7; Side 1 ]{{\includegraphics[width=6cm]{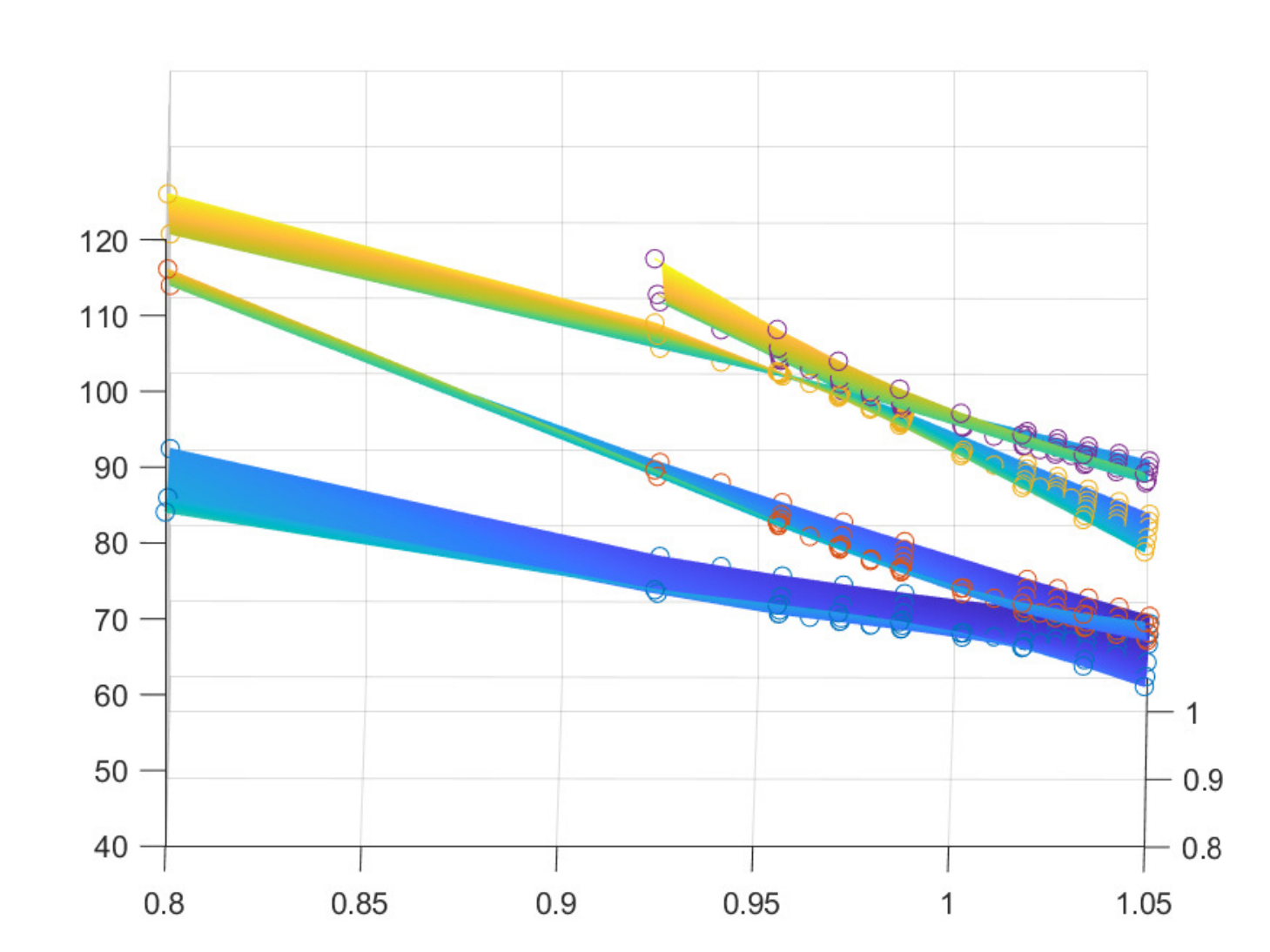}}}
    \qquad
    \subfloat[\centering Hypersurfaces 4,7; Side 2 ]{{\includegraphics[width=6cm]{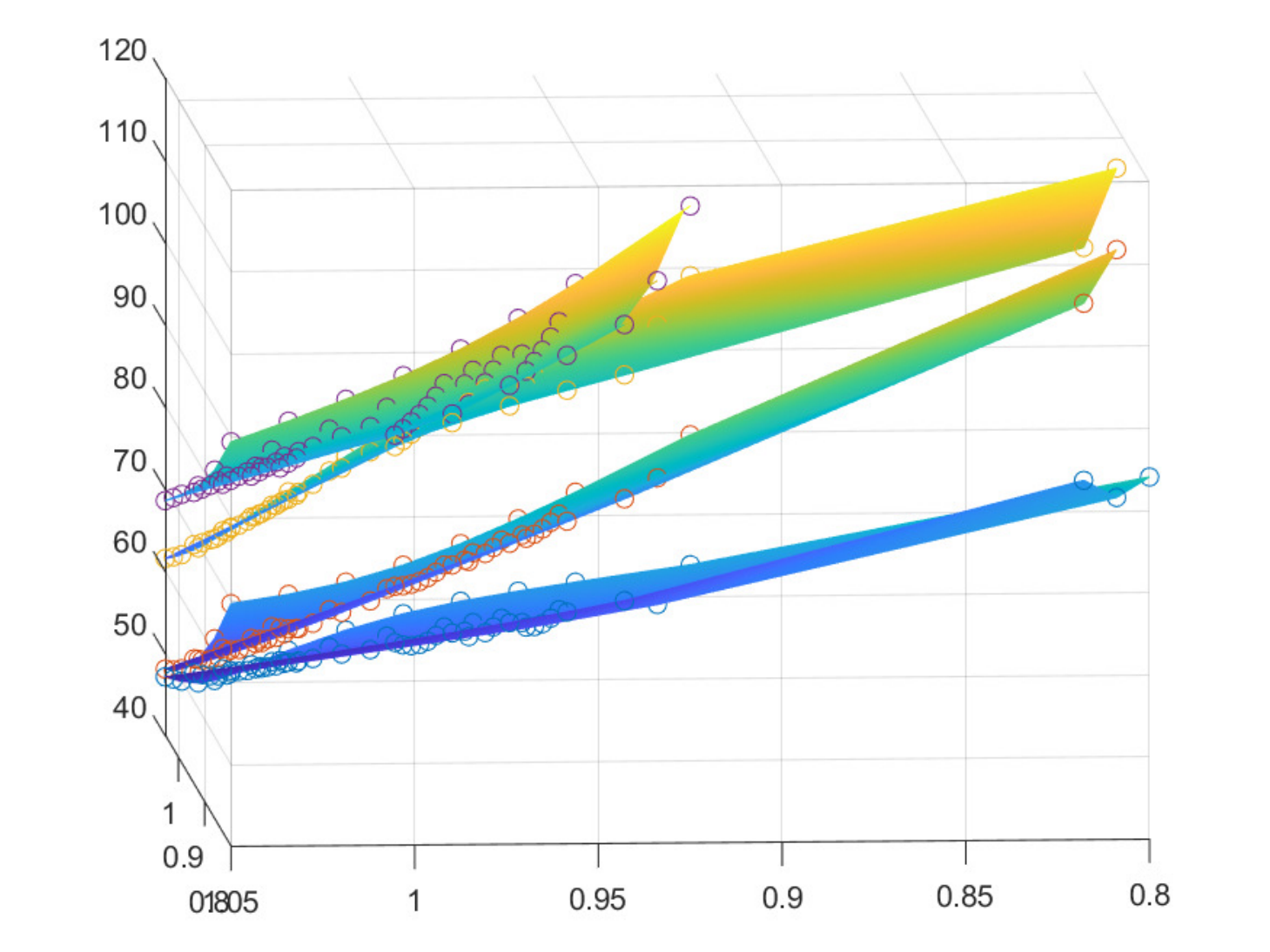}}}
    \qquad
    \subfloat[\centering Hypersurfaces 3,4 ]{{\includegraphics[width=6cm]{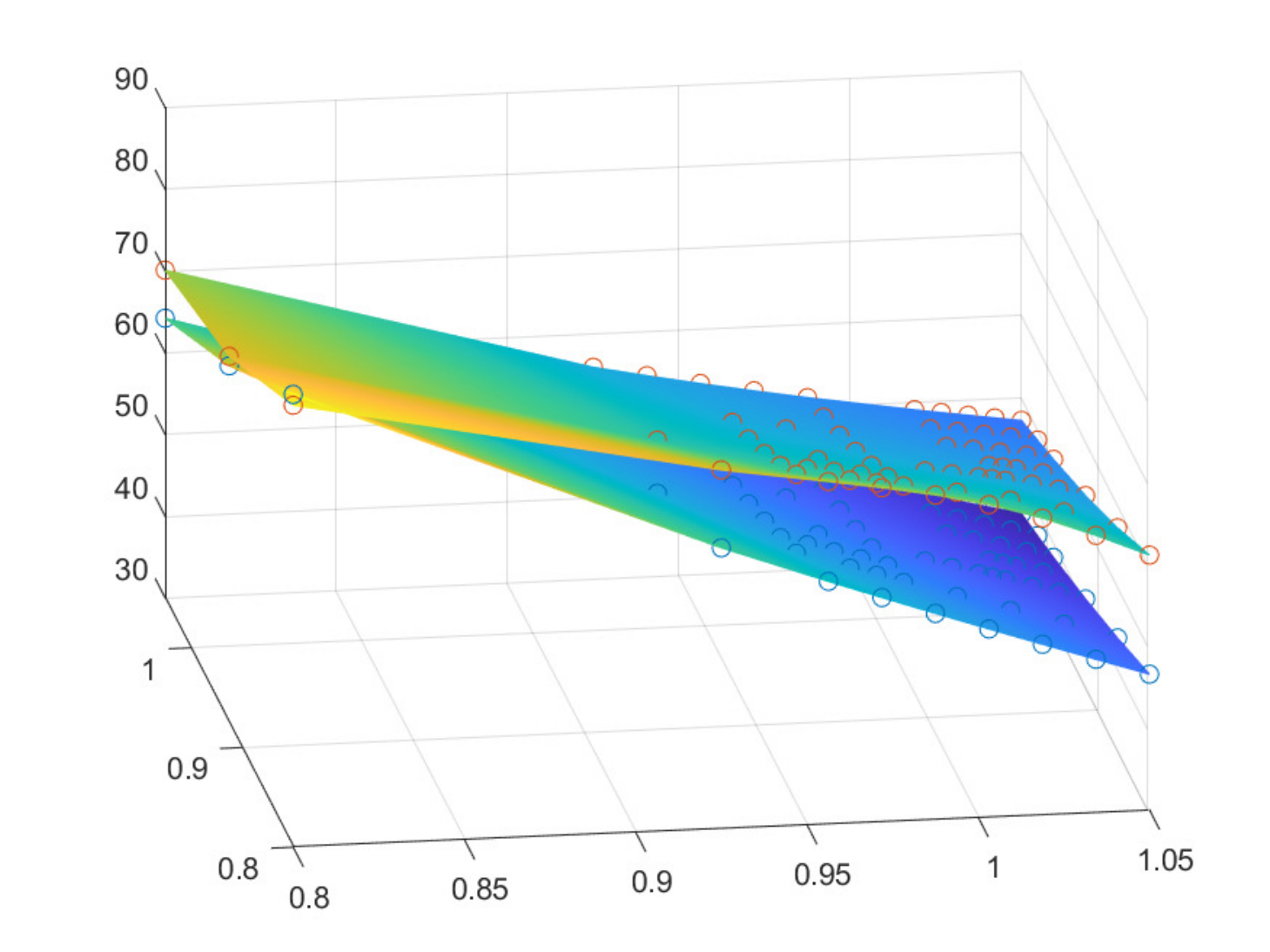}}}
    \caption{Final solution for the 2D example; Detailed view.}
    \label{fig:FinalSolution2DByBand}

\end{figure}

Given the initial grid $P^{(0)}$ being the $3\times3$ uniform lattice, and the tolerances $t_{\pi}=0.57$, $t_{\lambda}=0.015$, we let the adaptive algorithm run. The level-by-level output is displayed in Figure~\ref{fig:2D_grid_level_by_level}. The parametric grid produced by the adaptive algorithm is clearly non-uniform, and two regions in the parametric space can be recognized. At the left-half of $\mathcal M$, where no refinement happens, the eigenvalue hypersurfaces are well separated except for the crossing of the 3rd and 4th ones (see Figure~\ref{fig:FinalSolution2DByBand} (c)), which is however correctly identified by the a priori matching and subsequently certified by the a posteriori verification already at level 0. On the other hand, all the added points lie in the right-half of $\mathcal M$, even though no crossings are present. This refinement pattern is due to small gaps between the eigenvalue hypersurfaces (see Figure~\ref{fig:FinalSolution2DByBand} (a)-(b)). In particular, the extra grid points added at each level are meant to help the a posteriori test certify the identity a priori matching. 

\begin{table}
\begin{center}
\begin{tabular}{ |p{1cm}|p{2cm}|p{2.5cm}|p{2cm}|p{2.5cm}|   }
\hline
\multicolumn{5}{|c|}{Error By Level} \\
\hline
Level & Total no. of points & No. of wrongly matched points & No. of subintervals & No. of uncertified subintervals \\
\hline
0 & 9 & 7 & 12 & 4 \\
\hline
1  & 13 & 2 & 10 & 9 \\
\hline
2  & 22 & 1 & 22 & 11\\
\hline
3  & 36 & 4 & 36 & 10 \\
\hline
4  & 49 & 2 & 39 & 4 \\
\hline
5  & 59 & 0 & 34 & 2\\
\hline
6  & 64 & 0 & 20 & 0\\
\hline
\end{tabular}
\end{center}
\caption{Error table for the 2D example.}
\label{table:2DError}
\end{table}

Table~\ref{table:2DError} summarizes the level-by-level information. We can observe that, from level 5 on, five points are added to the grid, even though all the eigenpairs were correctly matched. These points serve the purpose of allowing the algorithm to fully certify all subintervals, and produce a correct solution after termination.
An important feature to notice is that some refinements can lead to an increase in the error (represented by the number of wrongly matched eigenpairs). This happens for example at level 3. Additional refinements at each level generate smaller subintervals. While increased accuracy is obtained, this doesn't translate into a continuous decay in the error. A correct solution is only guaranteed when all subintervals are certified. 

We conclude the section by bringing to the surface one last issue, namely how to propagate the matching information. In the one dimensional case, this issue was easily solved by defining the propagation direction from left to right. In the two dimensional case, the way to proceed is not clear anymore, and it becomes even more complicated for $d>2$.

Our technique relies on the creation of a path connecting all the points of the adapted grid. Such a path (which in principle is not guaranteed to be unique) must fulfill good properties (e.g., connected, no cycles are allowed) since the information must propagate distinctively. This step is performed in the code making use of the \textsc{Matlab} commands \verb|minspantree| and \verb|shortestpath|.

\begin{figure}[htpb]
\centering
\includegraphics[width=6cm]{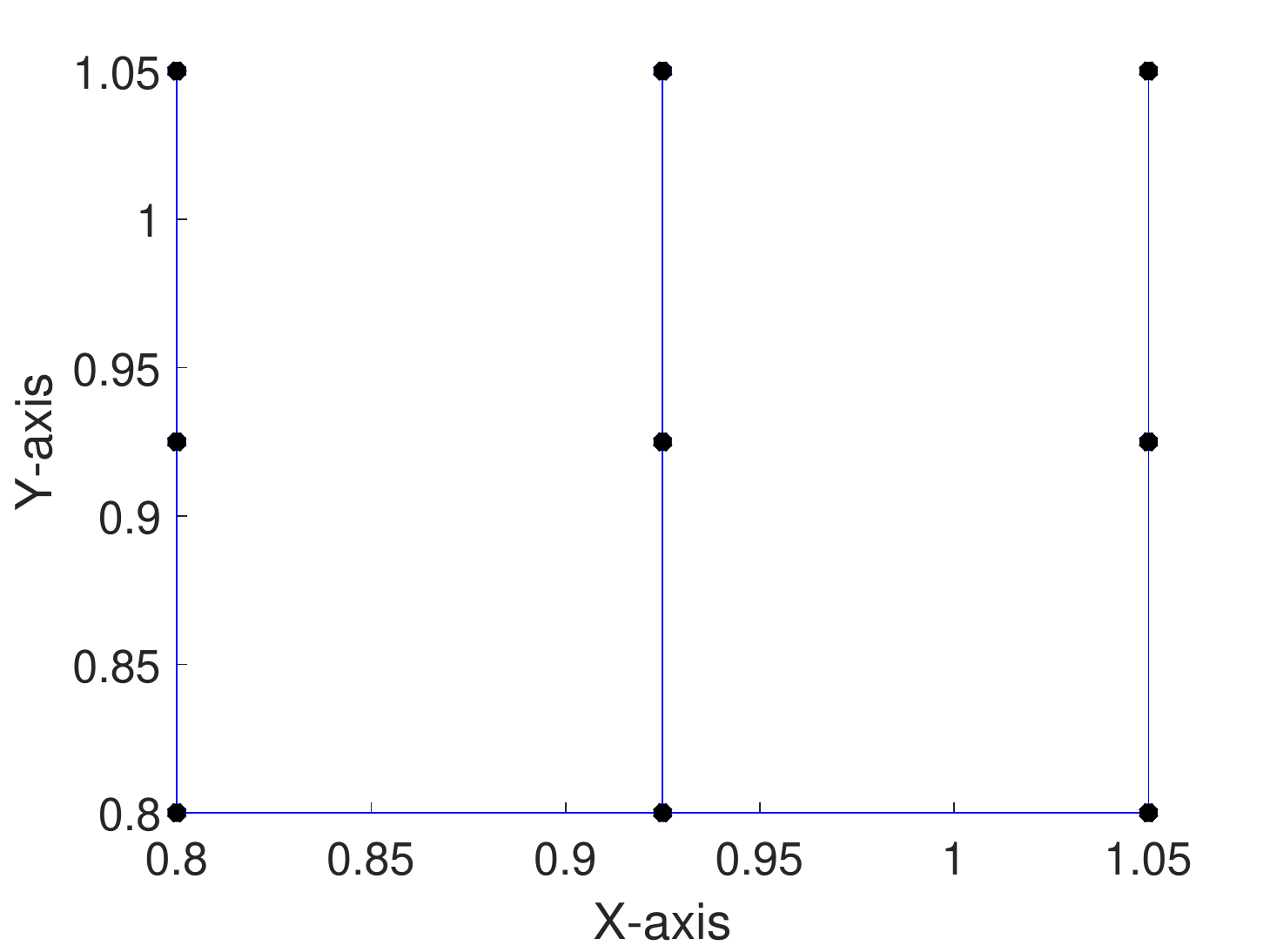}
\caption{Initial path.}
\label{fig:InitialPath}
\end{figure}

As an example, we display the initial path corresponding to level 0 in Figure~\ref{fig:InitialPath}. 
We generate this path by applying the \verb|minspantree| function on the graph 
$$G =\verb|graph|([1 \ 1 \ 2 \ 4 \ 2 \ 3 \ 5 \ 4 \ 5 \ 6 \ 7 \ 8],[2 \ 4 \ 5 \ 5 \ 3 \ 6 \ 6 \ 7 \ 8 \ 9 \ 8 \ 9]).$$ 
We then find our one-dimensional sub-paths using \verb|shortestpath|. In this case, the three sub-paths all start at (0.8, 0.8) and end at (0.8, 1.05), (0.925, 1.05) and (1.05, 1.05) respectively.

%%%%%%%%%%%%%%%%%%%%%%%%%%%%%%%%%%%%%%%%%%%%%%%%%%%%%%%%%%%%%%%%%%%%%%%%%%%%%%%%%%
\section{Conclusions}
\label{se:conclusions}
%%%%%%%%%%%%%%%%%%%%%%%%%%%%%%%%%%%%%%%%%%%%%%%%%%%%%%%%%%%%%%%%%%%%%%%%%%%%%%%%%%
The present paper introduces a novel adaptive algorithm for the numerical treatment of parametric eigenvalue problems arising from elliptic partial differential equations. It is composed of two phases: locally, we look at a specific subinterval and we decide (by means of an a priori matching followed by an a posteriori verification) whether to mark it for refinement or not; globally, we perform a sparse grid-based refinement step, which delivers an adapted grid in the parameter space refined where needed in order to detect the features and crossings of the eigenvalue hypersurfaces. Finally, a surrogate for the parameter-to-eigenvalue map is constructed simply by piecewise linear interpolation. Notably, the construction of a surrogate for the parameter-to-eigenfunction (or parameter-to-eigenvector) map is more delicate, and it is left for future investigations.

Even though the algorithm is written in an arbitrary dimension of the parameter space, numerical examples are performed in 1D and 2D, only, with the scope of attesting the validity and verifying the performances of the proposed numerical scheme. Higher-dimensional numerical tests will be presented in a forthcoming contribution.

This work paves the way towards the treatment of stochastic eigenvalue problems, i.e., eigenvalue problems arising from elliptic partial differential equations with random coefficients. In recent years a huge effort has been made in the study of uncertainty quantification (UQ) techniques for the source problem, and various methods have been developed (we mention, e.g., the Monte Carlo method~\cite{RobertCasella}, non-intrusive and Galerkin methods~\cite{LeMaitreKnio2010} and perturbation methods~\cite{bn_20,bnk_16,bn_14}. However, the field of stochastic eigenvalue problems is still quite unexplored, and we believe that the combined use of UQ techniques together with the algorithm proposed here represents a promising way to go.

\section*{Acknowledgements}
The work of D.\ Boffi was supported by the Competitive Research Grants Program CRG2020 ``Synthetic data-driven model reduction methods for modal analysis'' awarded by the King Abdullah University of Science and Technology (KAUST). D.\ Boffi is member of the INdAM Research group GNCS and his research is partially supported by IMATI/CNR and by PRIN/MIUR.
F.\ Bonizzoni is member of the INdAM Research group GNCS and her work is part of a project that has received funding from the European Research Council ERC under the European Union's Horizon 2020 research and innovation program (Grant agreement No.~865751).

\end{document}